\documentclass {article}

\usepackage{pinlabel}

\usepackage{graphicx}
\usepackage{amsthm}

\newcommand{\pic}[2]{\raisebox{-.5\height}{\includegraphics[scale=#2]{#1}}}

\def\kernel{\pic{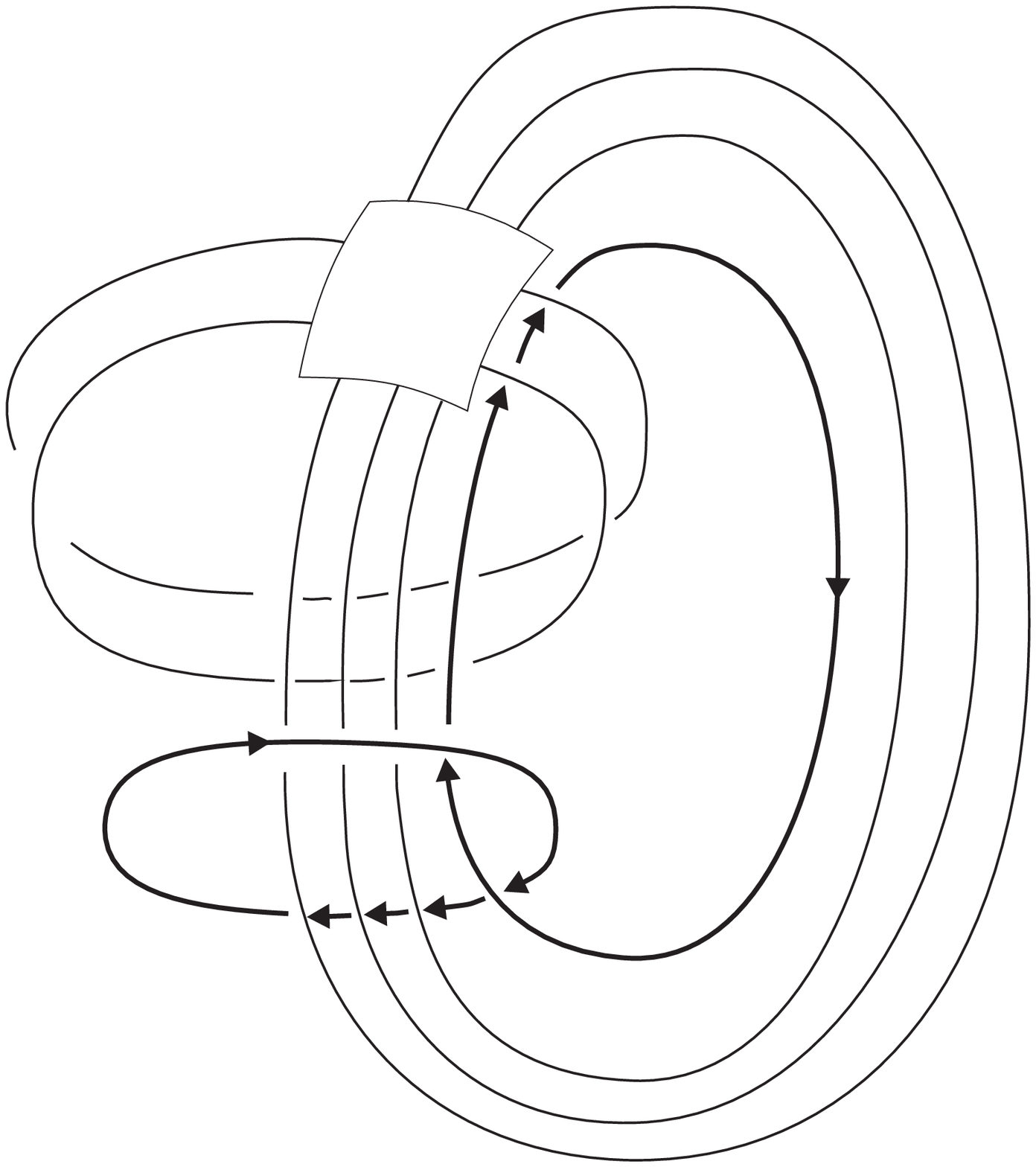}{.250}}

\def\plainweavecell{\pic{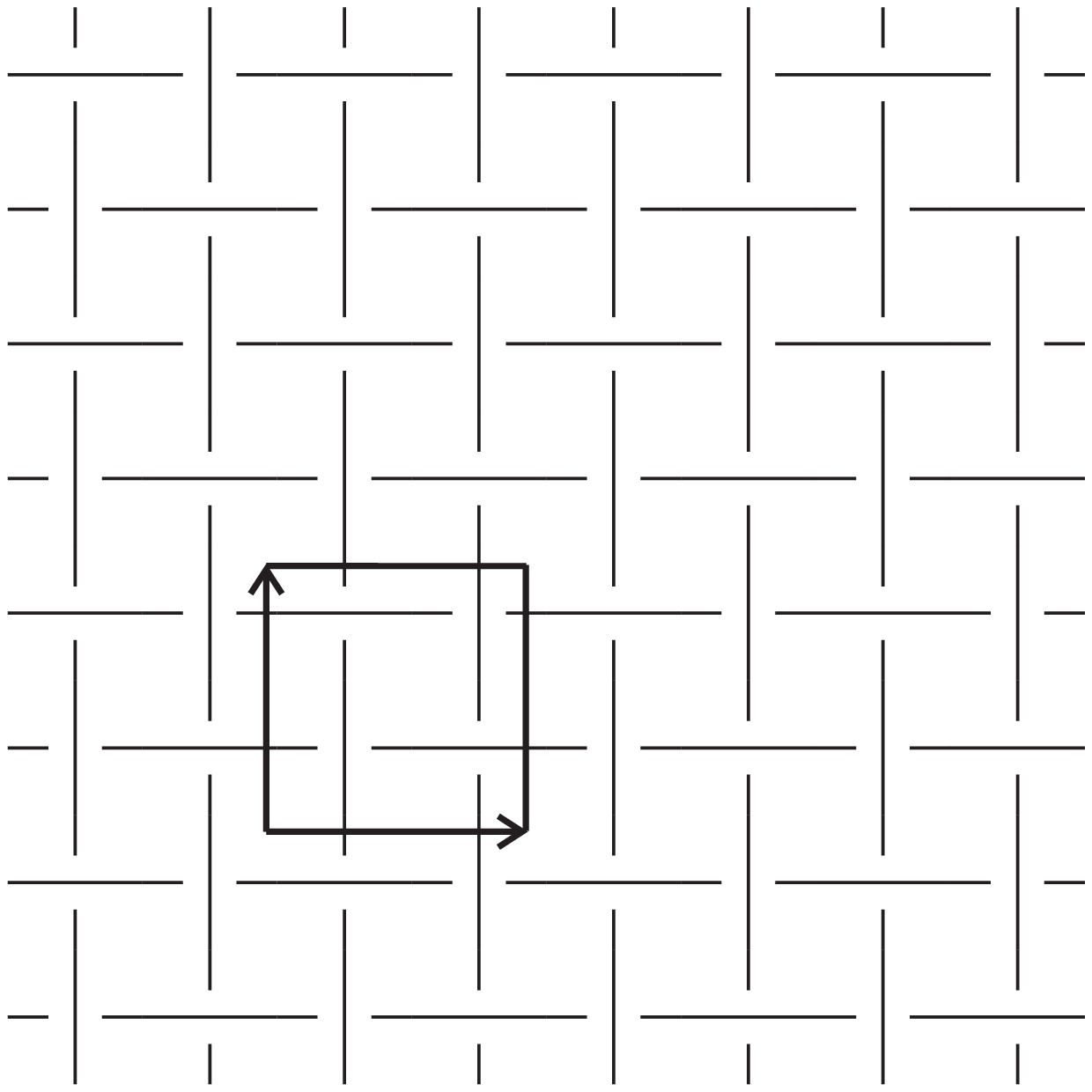}{.30}}
\def\plainweavekernel{\pic{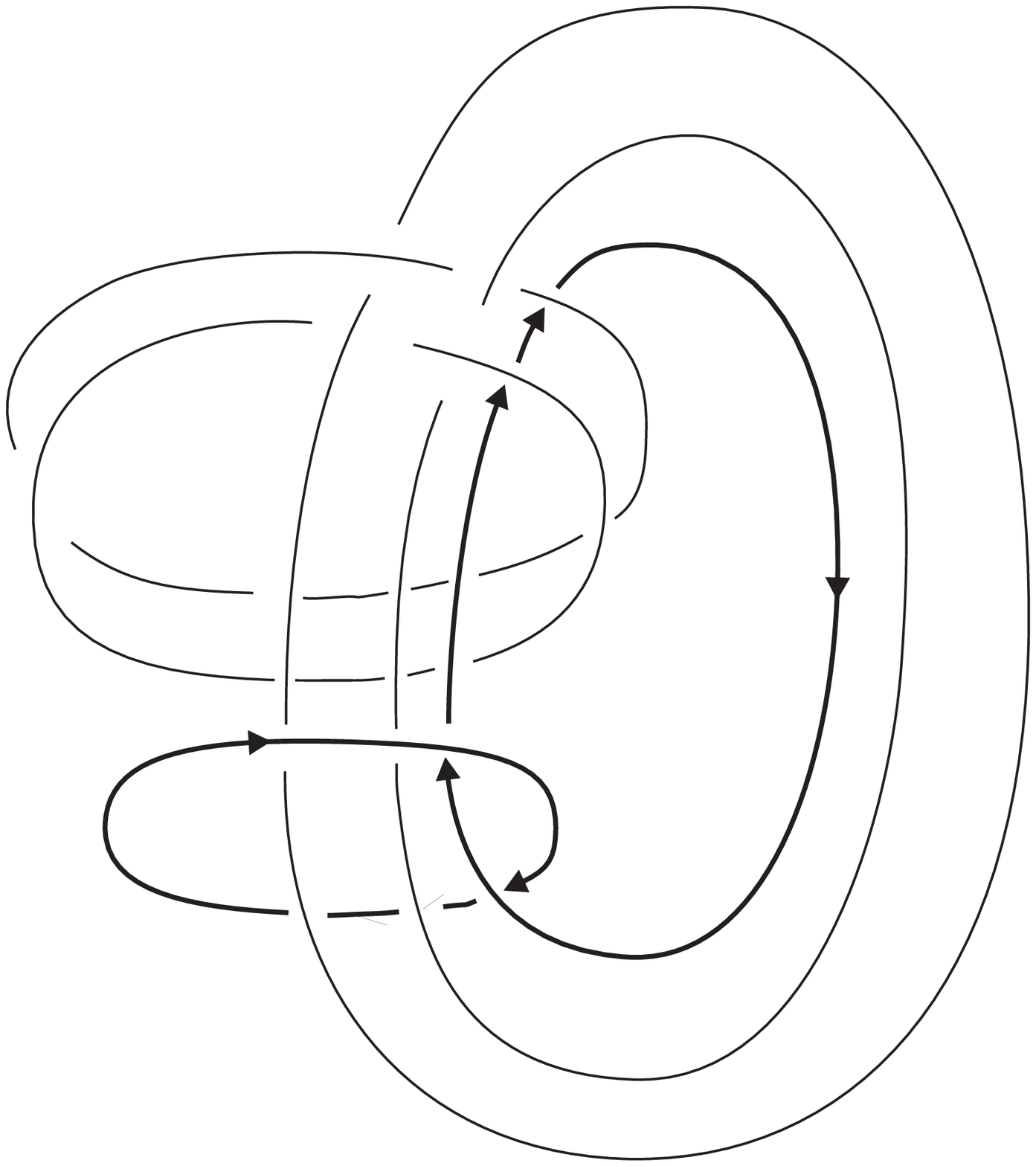}{.220}}

\def\lenoweavecell{\pic{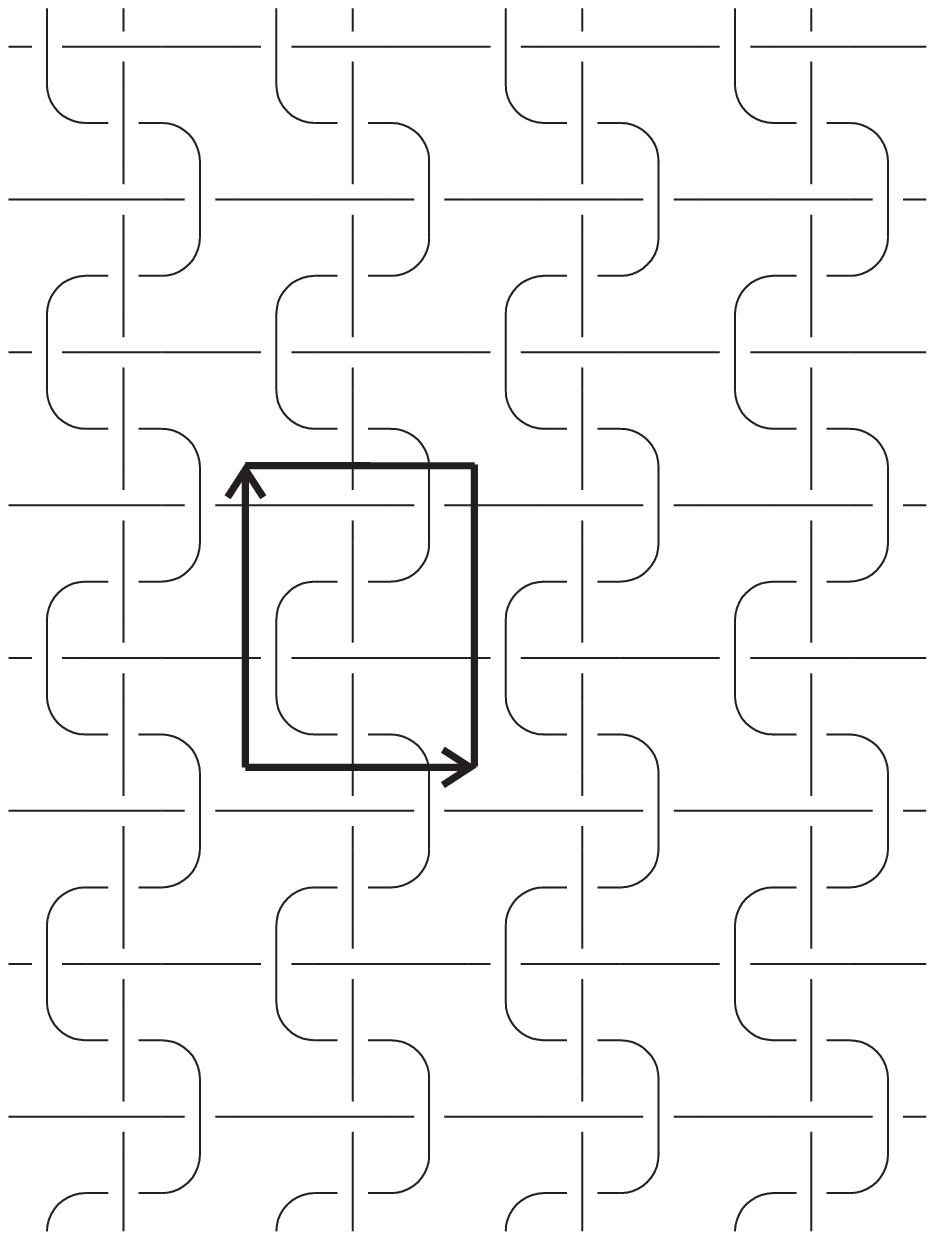}{.40}}
\def\lenoweavekernel{\pic{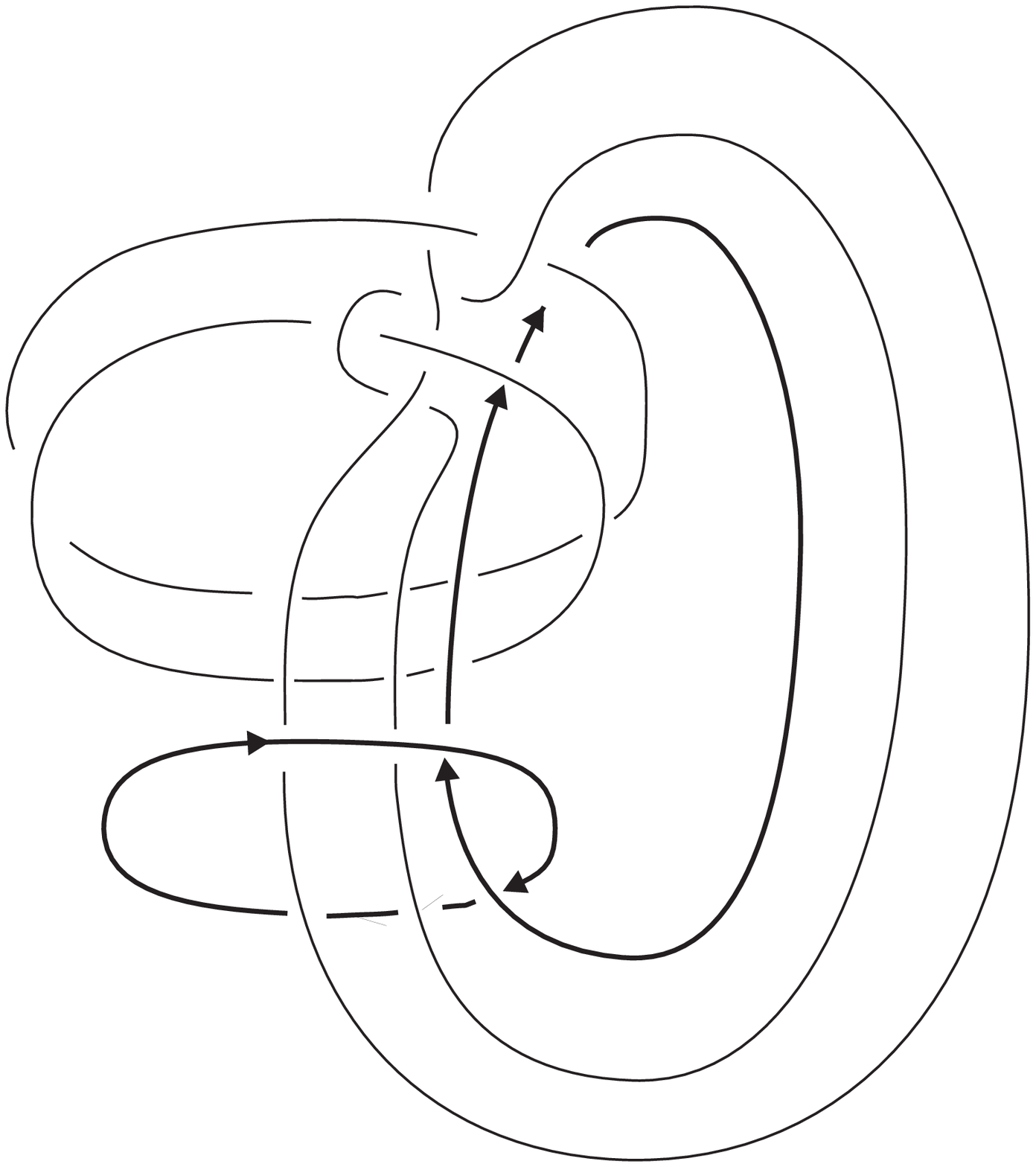}{.220}}

\def\fouraxialcell{\pic{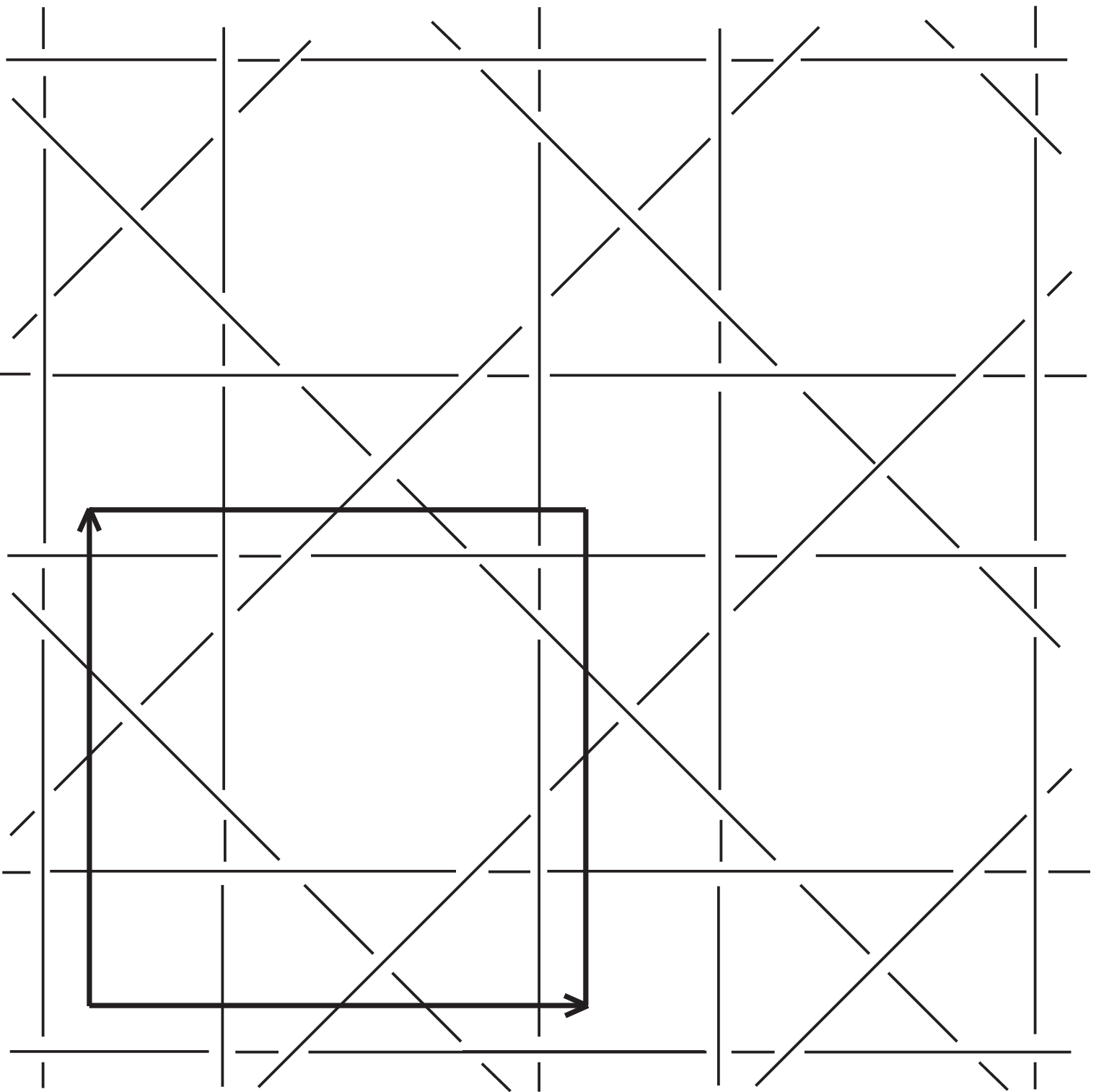}{.30}}
\def\quatroaxialkernel{\pic{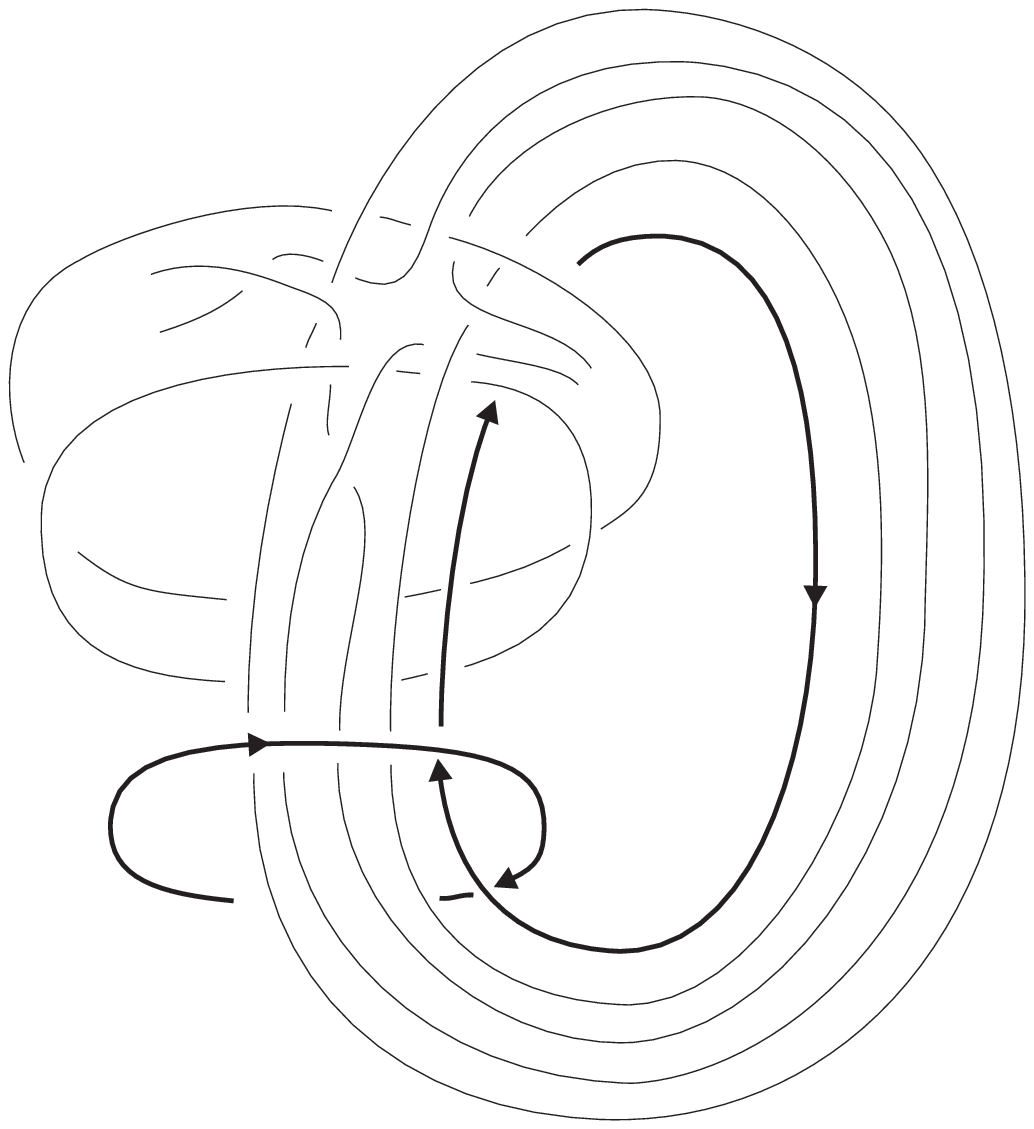}{.4}}

\def\twillcell{\pic{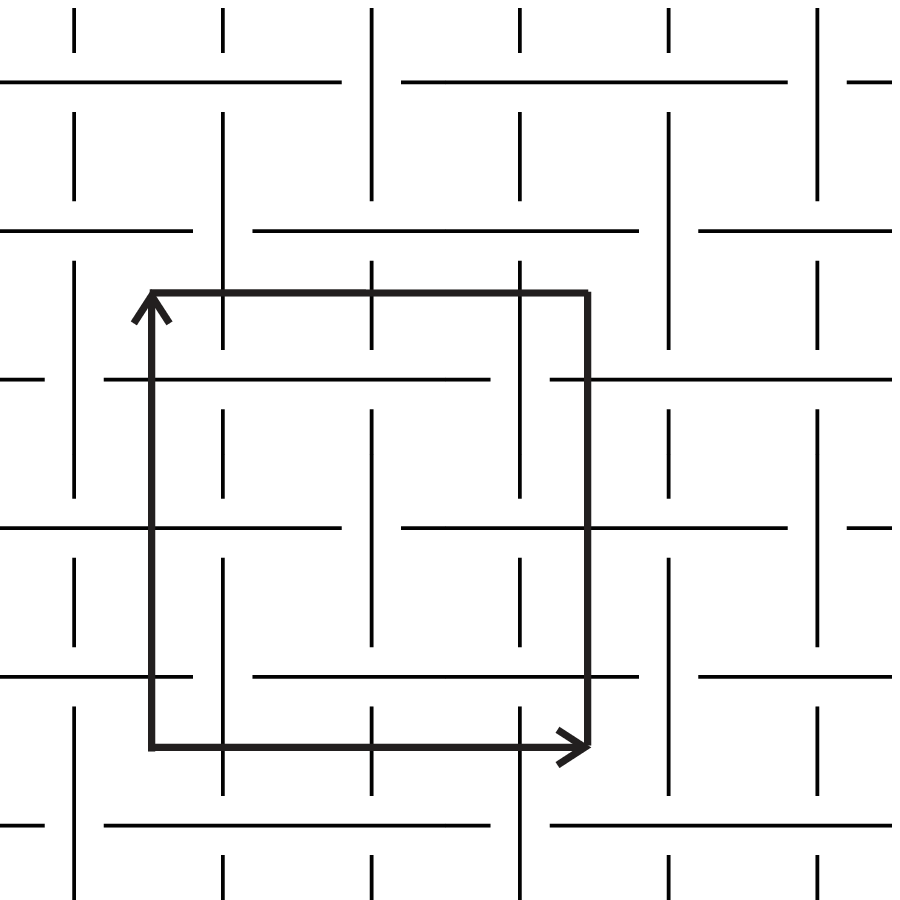}{.40}}
\def\twillmincell{\pic{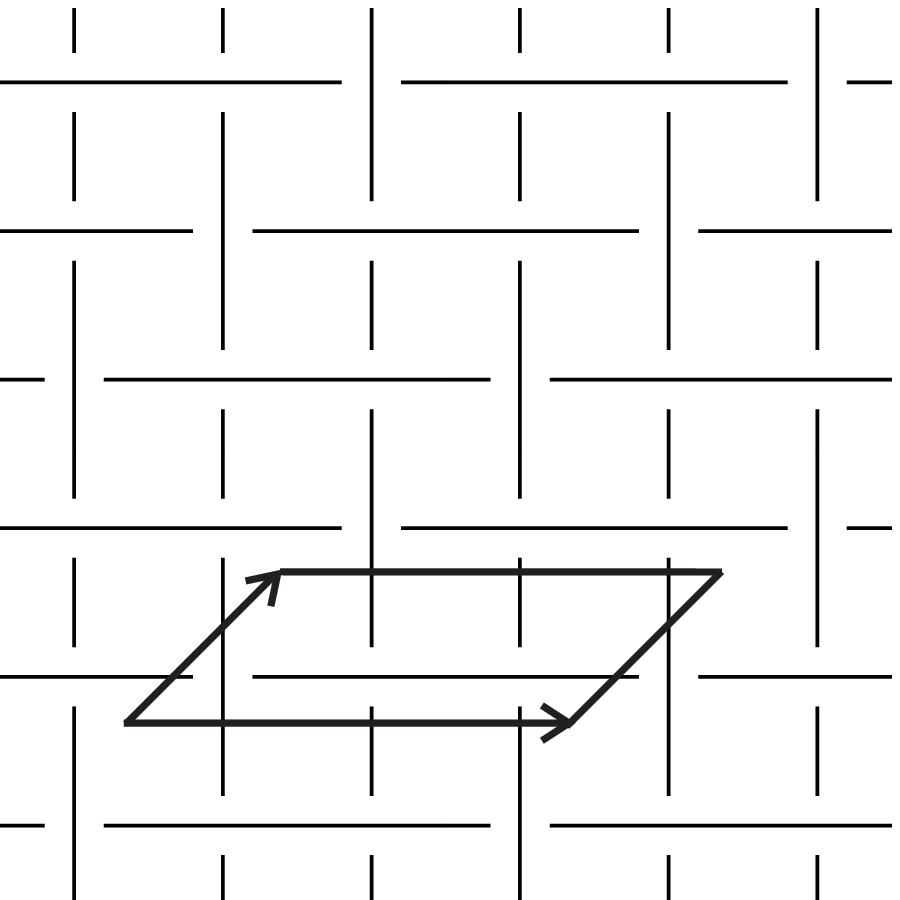}{.40}}
\def\singlejerseycell{\pic{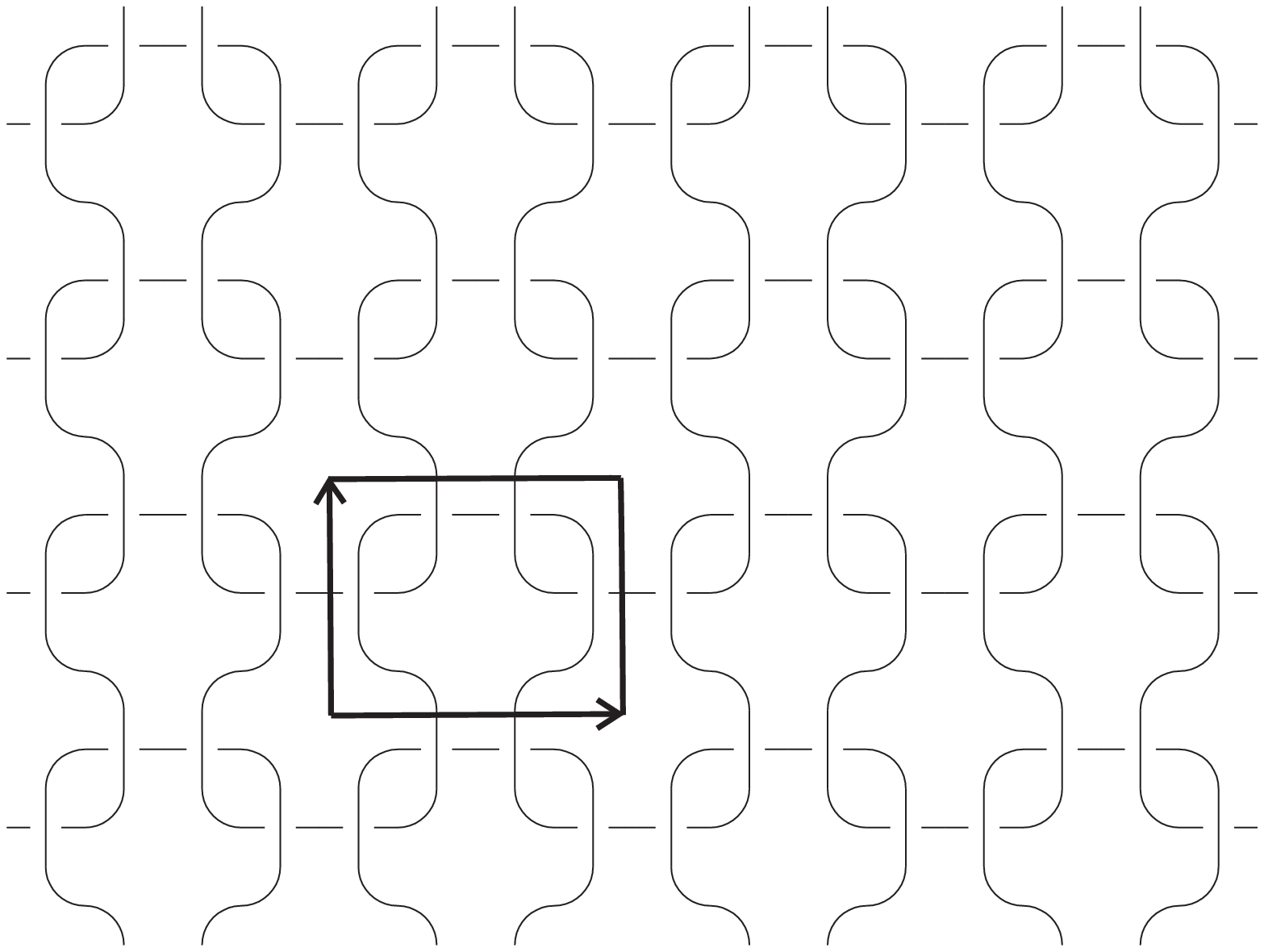}{.30}}
\def\singlejerseykernel{\pic{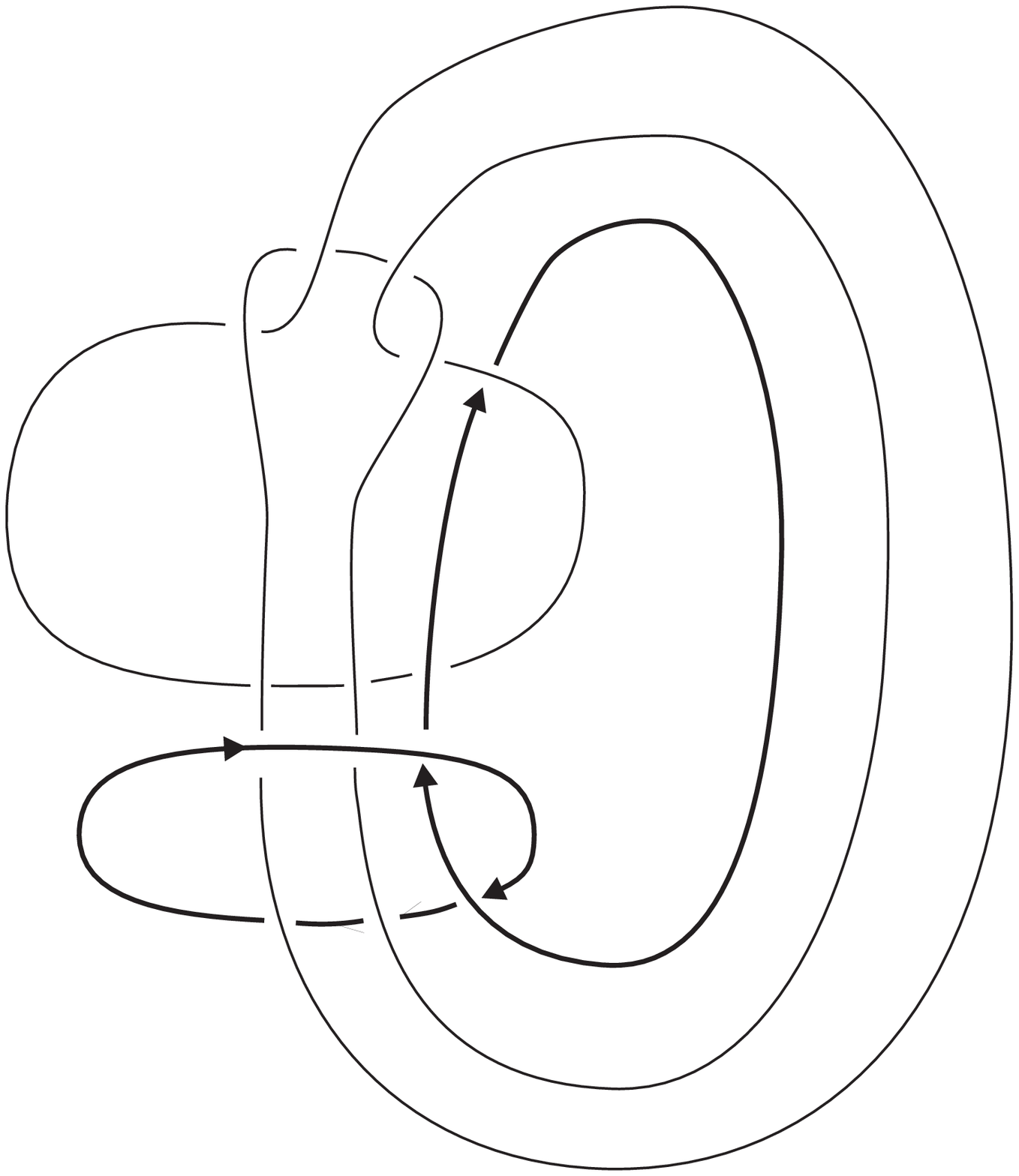}{.220}}

\def\singlejerseycellskew{\pic{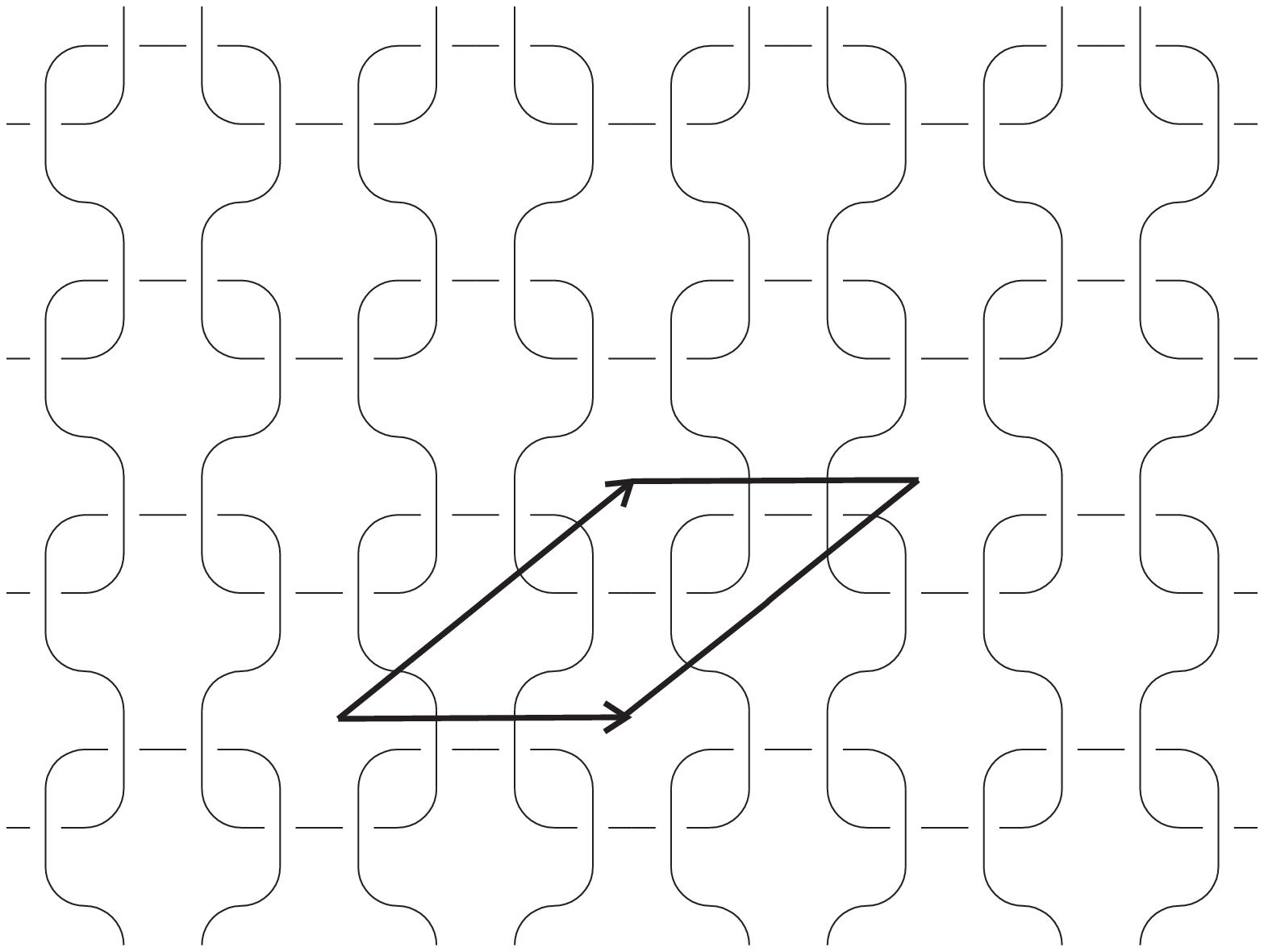}{.30}}
\def\singlejerseykernelskew{\pic{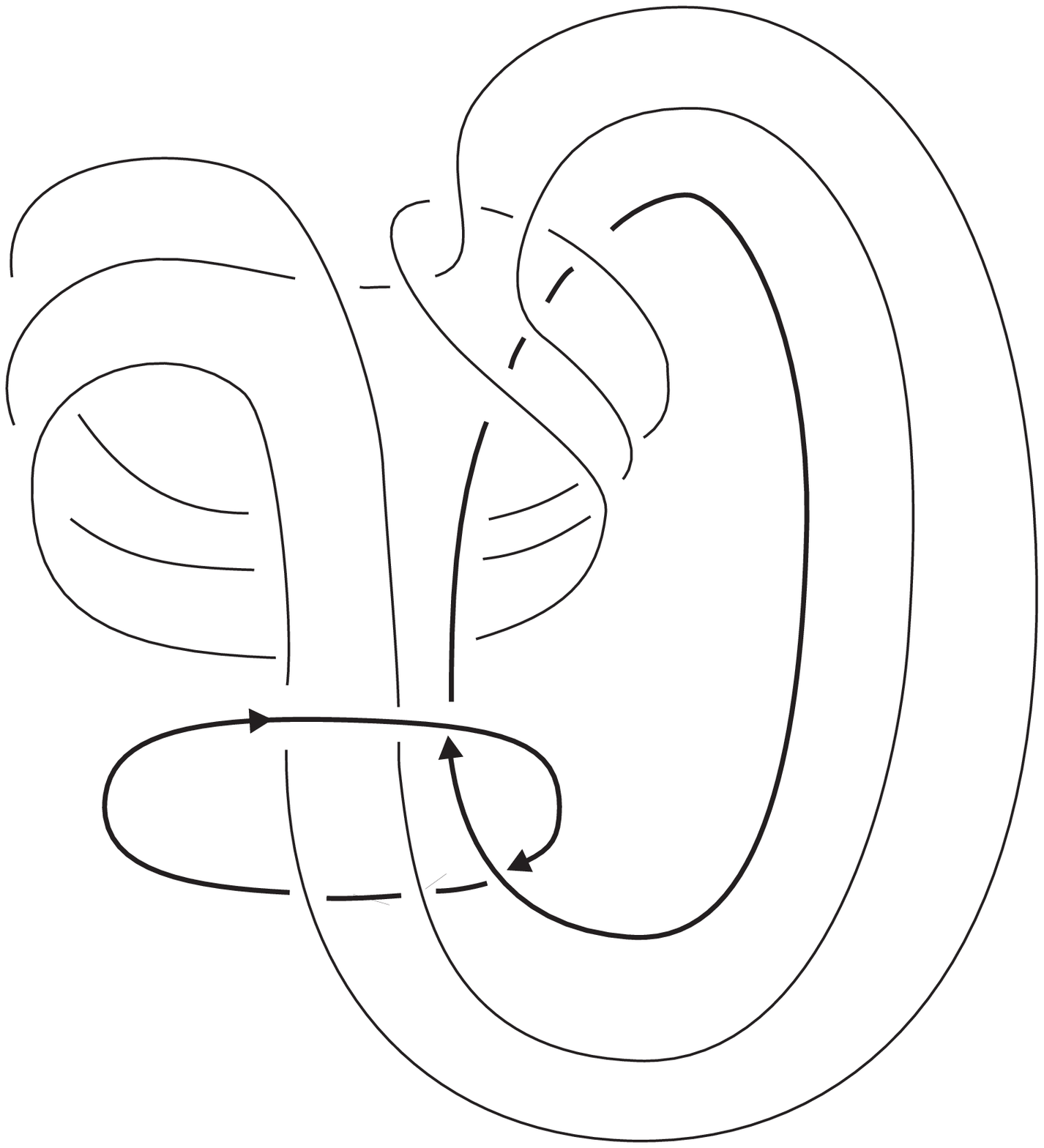}{.220}}

\def\singlejerseywarpweft{\pic{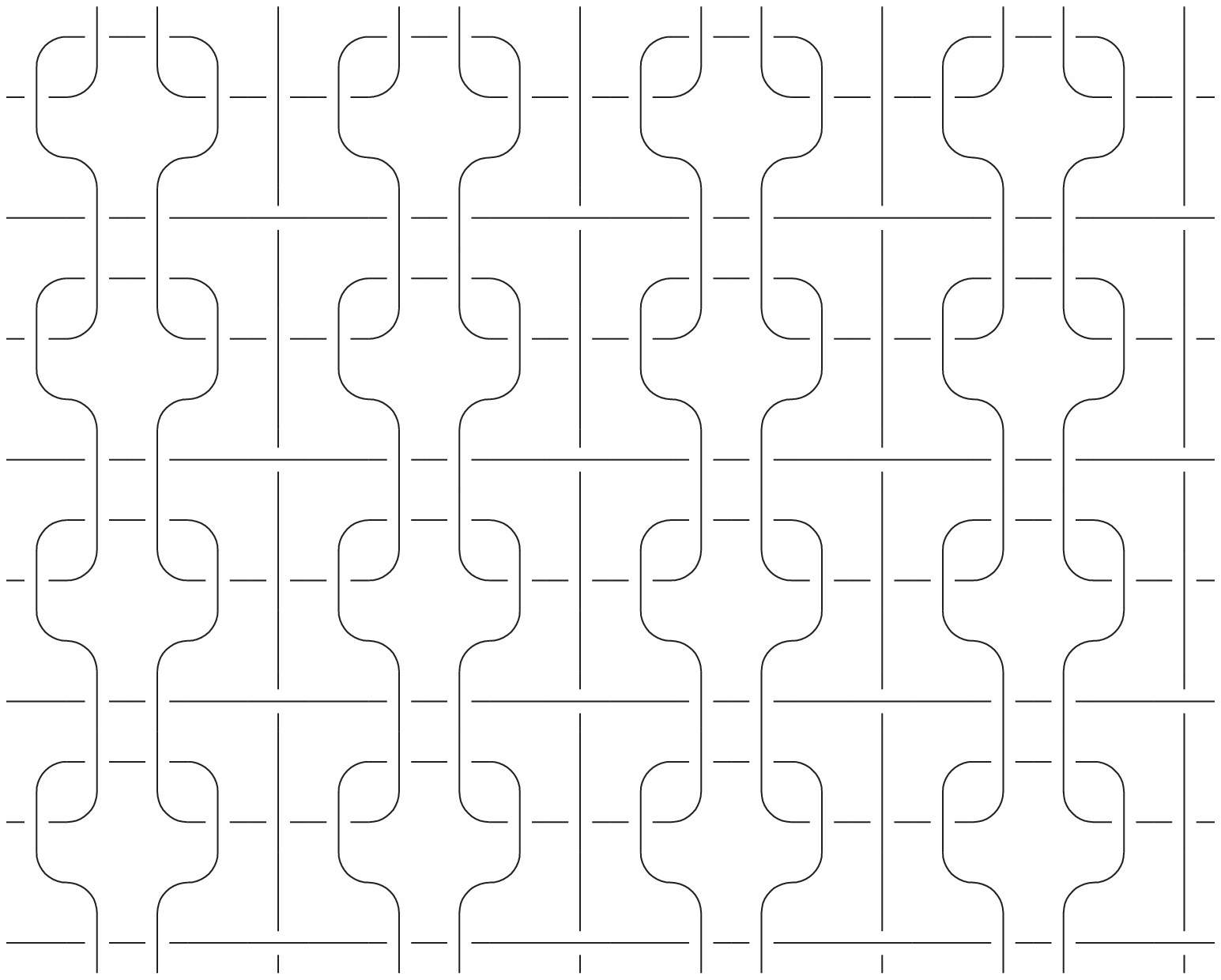}{.30}}

\def\triaxial{\pic{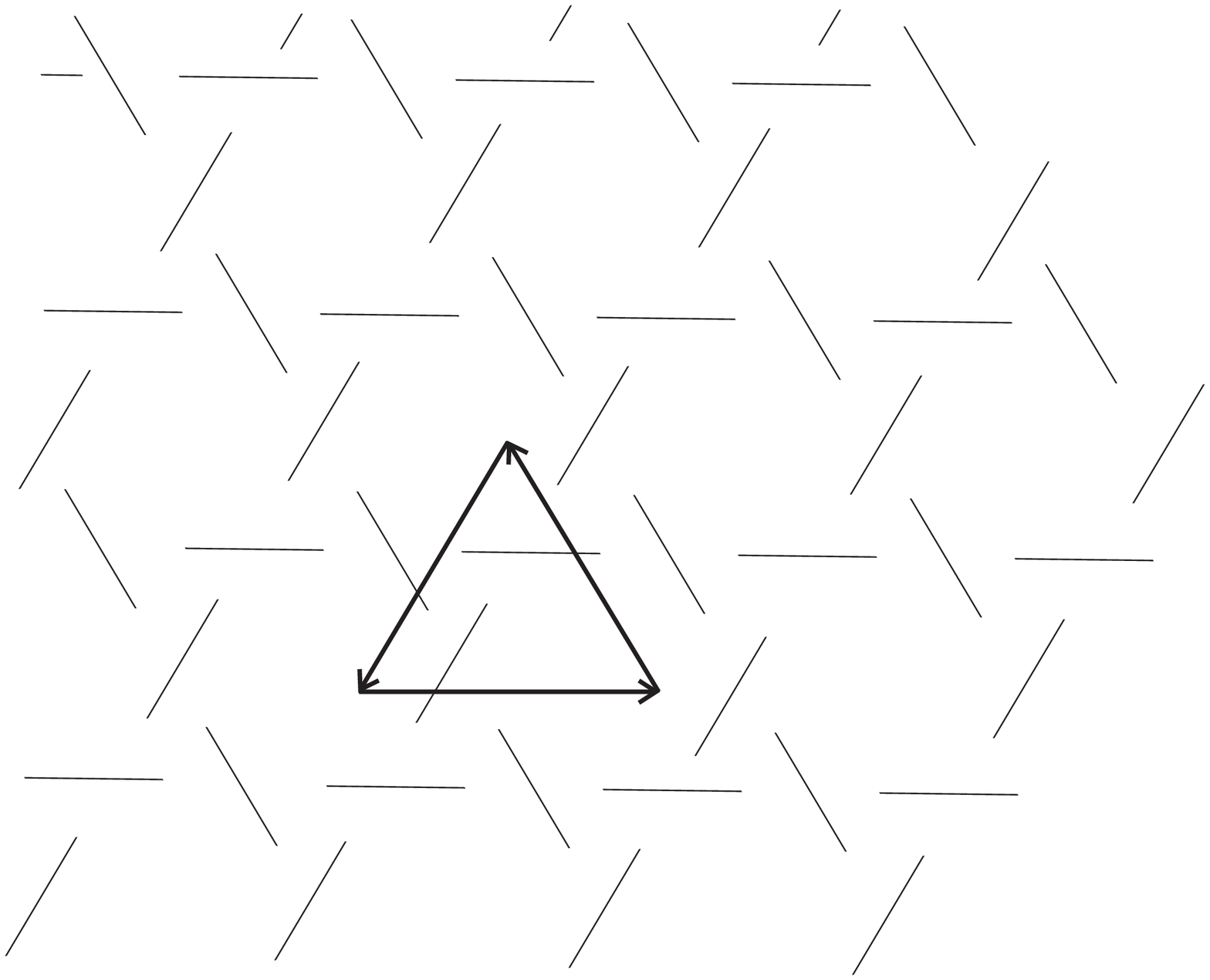}{.30}}

\def\triaxialdata{\pic{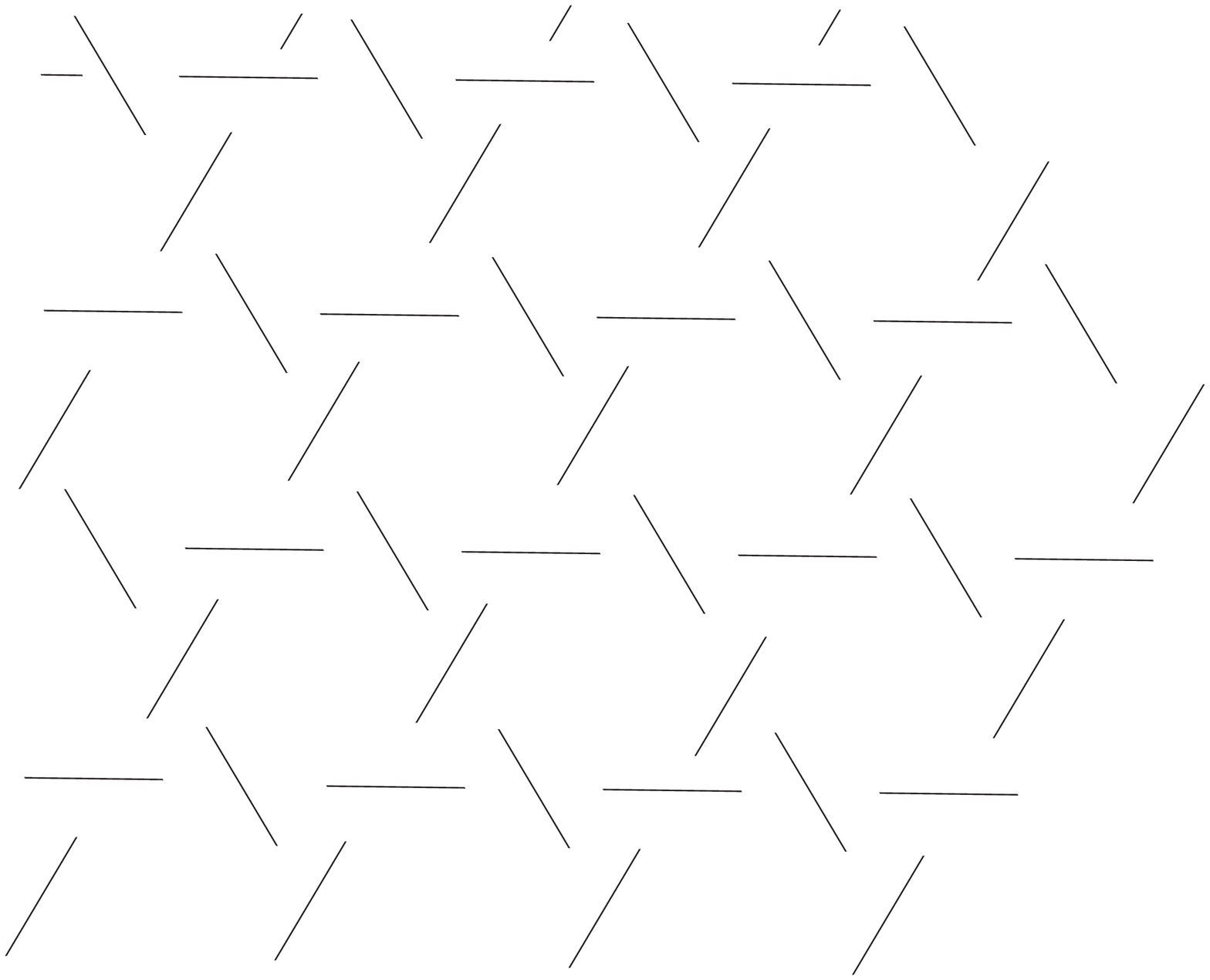}{.30}}

\def\singlejerseyback{\pic{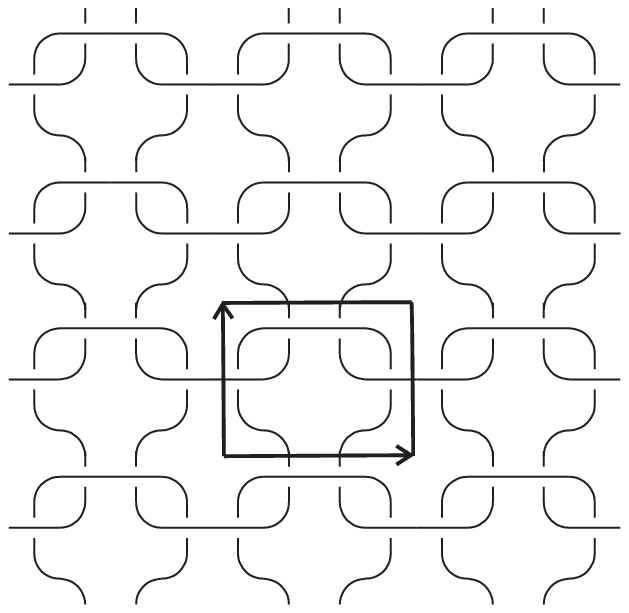}{.50}}

\def\fishnetcell{\pic{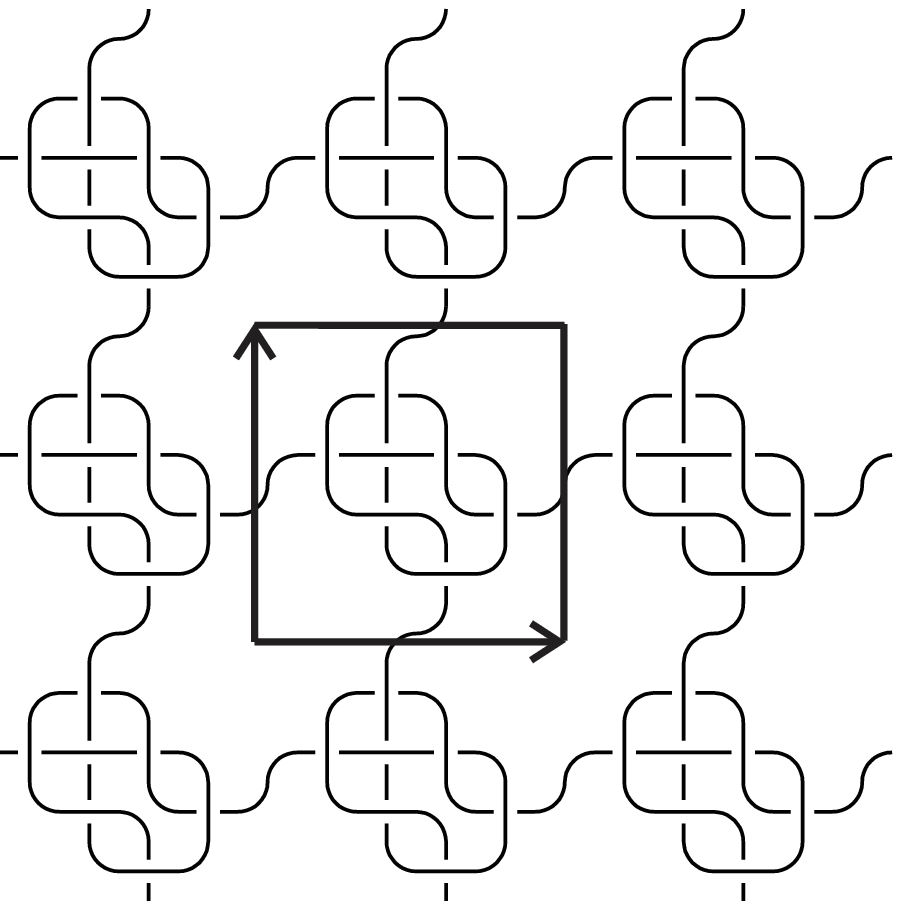}{.40}}

\def\chaincell{\pic{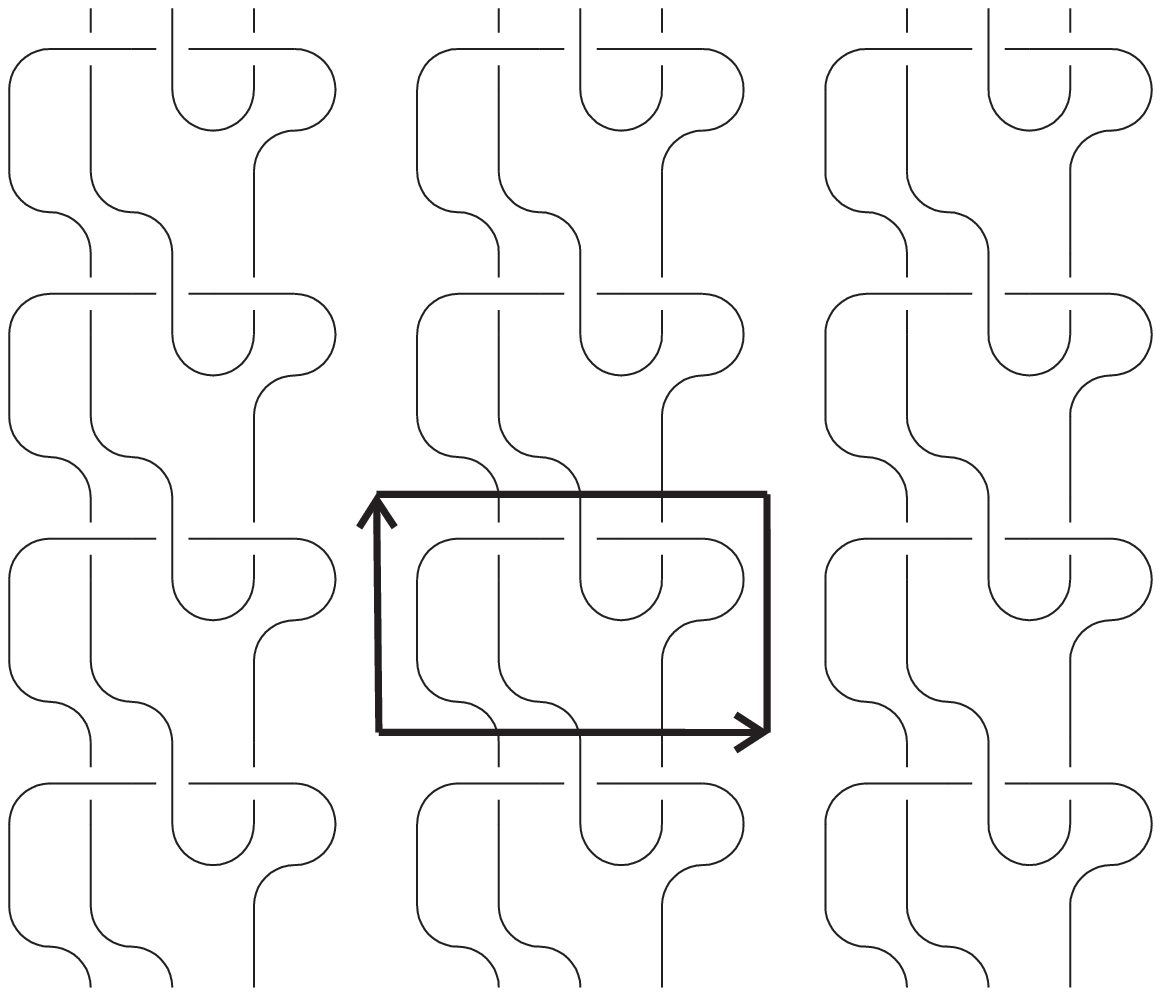}{.40}}

\def\chainmail{\pic{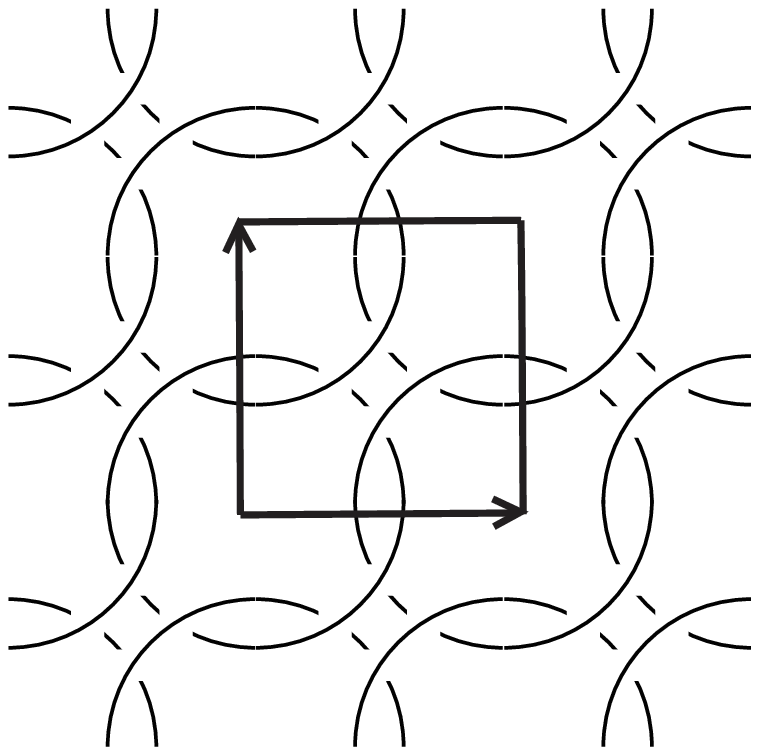}{.40}}

\def\fabric{\pic{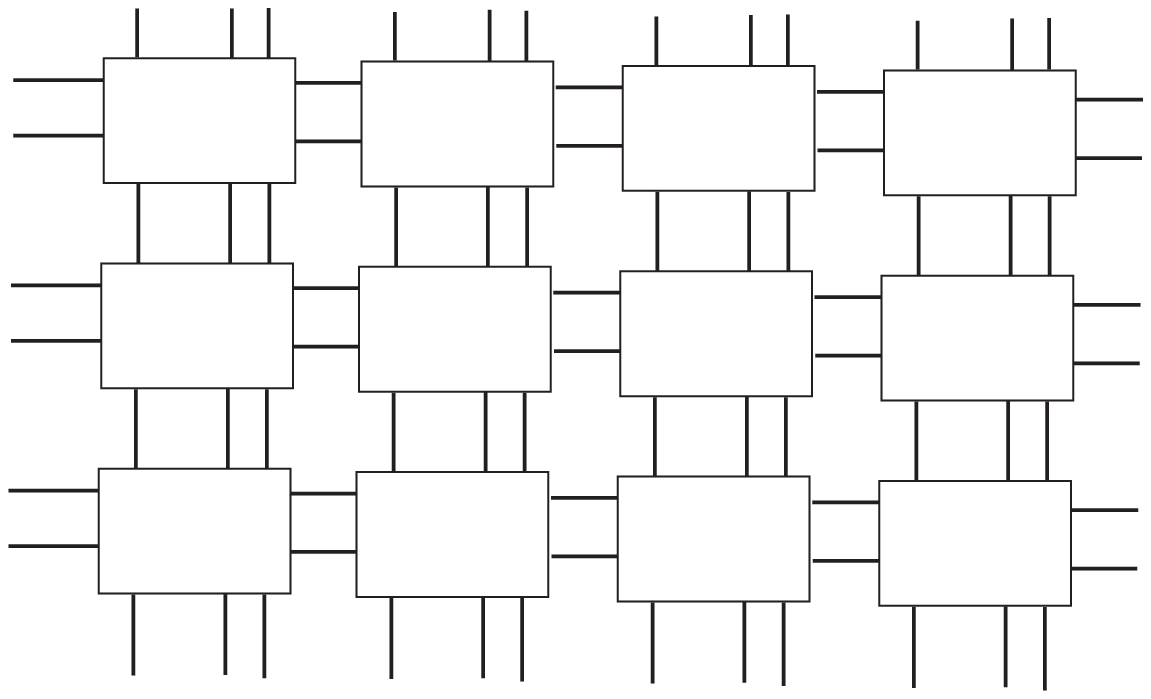}{.45}}

\newcommand{\bc}{\begin{center}}
\newcommand{\ec}{\end{center}}

\newcommand{\be}{\begin{equation}}
\newcommand{\ee}{\end{equation}}

\newcommand{\beqn}{\begin{eqnarray*}}
\newcommand{\eeqn}{\end{eqnarray*}}

\newcommand{\ds}{\displaystyle}
\newcommand{\x}{\times}
\newcommand{\ha}{\frac{1}{2}}

\newtheorem{theorem}{Theorem}
\newtheorem{corollary}[theorem]{Corollary}
\newtheorem{lemma}[theorem]{Lemma}

\newenvironment{definition}{\par\smallskip%
\noindent\textbf{Definition.}\  }%
{\par\smallskip}

\theoremstyle{remark}
\newtheorem{example}{Example}
 \newtheorem{remark}{Remark}

\begin{document}

\bc{\Large\bf Doubly periodic textile patterns\\[3mm]}
{\sc H. R. Morton {\rm and } S. Grishanov\\[2mm]}
 {\small \sl Department of Mathematical Sciences, University of Liverpool,\\ Peach Street, Liverpool L69 7ZL.\\[2mm] 
Textile Engineering and Materials Research Group\\ De Montfort University, The Gateway, Leicester LE1 9BH.}

\ec

\begin{abstract} Knitted and woven textile structures are examples of doubly periodic structures in a thickened plane made out of intertwining strands of yarn. Factoring out the group of translation symmetries of such a structure gives rise to a link diagram in a thickened torus, as in \cite{Grishanov}. Such a diagram on a standard torus in $S^3$ is converted into a classical link by including two auxiliary components which form the cores of the complementary solid tori. The resulting link, called a \emph{kernel} for the structure, is determined by a choice of generators $u,v$ for the group of symmetries.

A normalised form of the multi-variable Alexander polynomial of a kernel is used to provide polynomial invariants of the original structure which are essentially independent of the choice of generators $u$ and $v$.
It gives immediate information about the existence of closed curves and other topological features in the original textile structure. Because of its natural algebraic properties under coverings we can recover the  polynomial for kernels based on a proper subgroup  from the polynomial derived from the full symmetry group of the structure.  This enables  two structures to be compared at similar scales, even when one has a much smaller minimal repeating cell than the other.

 Examples of simple traditional structures are given, and their Alexander data polynomials are presented to illustrate the techniques and results.
\end{abstract}

\section{Introduction}

Textiles represent a diverse class of commonly used materials with specific structural
properties which impose a number of restrictions on the mutual position of constituting threads
and the way in which they are intermeshed. Unlike general knots and links, textile structures
cannot contain closed components or knots tied on their threads; they must be structurally
coherent, i.e. fabrics cannot contain non-interlaced threads, disconnected strips or layers.

The development of specific tools to identify such forbidden elements in the fabric structure forms the focus of this paper.

\subsection{Representing a fabric by a link}
In \cite{Grishanov} Grishanov, Meshkov and Omelchenko introduced the idea of representing a
fabric with a repeating (doubly periodic) pattern by a knot  diagram on a torus, having made a
choice of a unit cell for the repeat of the  pattern. Algebraic invariants of this diagram
based on the Jones polynomial were used to associate  a polynomial to the fabric which was
independent of the choice of unit cell, so long as a minimal choice of repeating cell was made.
In this paper we enhance the nature of the diagram used to represent the fabric by including
two further auxiliary curves, to produce a link in the $3$-dimensional sphere $S^3$ from which the original fabric can
be recovered. We use the multivariable Alexander polynomial of the resulting link so as to
strengthen the information available about topological properties of the fabric, and remove the
need to work with a minimal choice of repeating cell.

We use the term \emph{fabric}  to mean a doubly periodic oriented plane knot diagram,
consisting of coloured strands with at worst simple double point crossings, up to the classic
Reidemeister moves.  A fabric gives rise to  a link diagram on the   torus $T^2\cong S^1\x S^1$  by
choosing a repeating cell in the pattern and splicing together the strands where they cross
corresponding edges to form the diagram on the torus.  A link in $S^3$ with two further
auxiliary components $X$ and $Y$ is constructed by placing the torus in $S^3$ as a standard
torus and including the core curves on each side of the torus in addition to the curves forming
the diagram on the torus.  We make the convention that the  curve $X$ lies on the side of the
torus towards the \emph{face} of the original fabric, and the curve $Y$ lies on the side of the
torus towards the \emph{back} of the fabric. The resulting link, with the distinguished choice
of curves $X$ and $Y$, will be called a \emph{kernel} for the fabric.

 We assume that our fabric
 lies in a thickened plane, and is invariant under a discrete group $G$ generated by two
 independent translations. The quotient of the thickened plane by the action of $G$ is then
 a thickened torus $T^2\x I$, bounded by two tori corresponding to the face and the back of
 the fabric. In forming a kernel of the fabric we have made a choice of embedding of this
 thickened torus in $S^3$, determined by an explicit choice of two generators $u$ and $v$
 for $G$ along the edges of our chosen unit cell.

\begin{figure}[ht!]
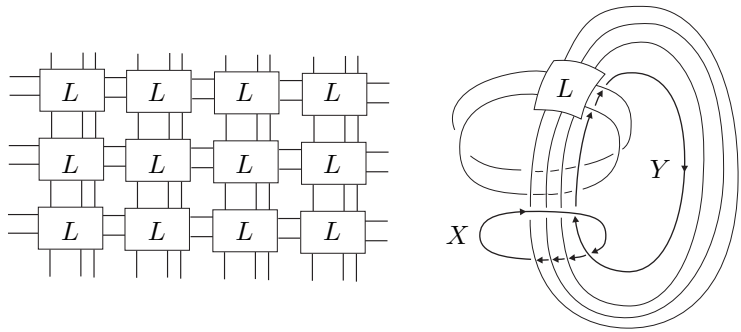

\bc
{\labellist
\pinlabel {$L$} at 122 725
\pinlabel {$L$} at 196 725
\pinlabel {$L$} at 271 725
\pinlabel {$L$} at 346 725
\pinlabel {$L$} at 122 664
\pinlabel {$L$} at 196 664
\pinlabel {$L$} at 271 664
\pinlabel {$L$} at 346 664
\pinlabel {$L$} at 122 606
\pinlabel {$L$} at 196 606
\pinlabel {$L$} at 271 606
\pinlabel {$L$} at 346 606
\endlabellist\fabric} \qquad
{\labellist
\pinlabel {$L$} at 176 504
\pinlabel {$X$} at 15 282
\pinlabel {$Y$} at 320 382
\endlabellist\kernel}\ec
\caption{A kernel for a fabric} \label{fabric}
\end{figure}

A schematic view of a fabric, with its face uppermost, and a kernel for it, are shown in Figure
\ref{fabric}.

\subsection{Some examples}
Some traditional fabrics are shown in the following figures, with a choice of repeating cell
indicated, and the resulting kernel in each case. A wide variety of fabrics can be found in the
books by Watson, \cite{Watson1,Watson2},  and Spencer, \cite{Spencer01}, and some of the primary
structural elements are described in \cite{Emery94}.

\begin{figure}[ht!]
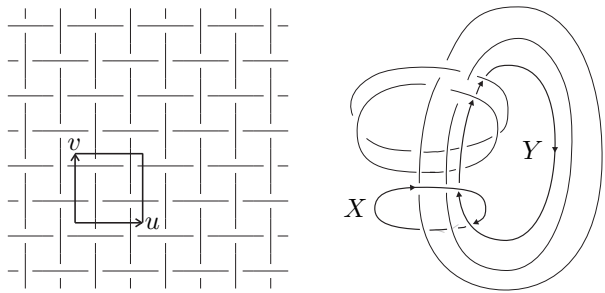

\bc
{\labellist
\pinlabel {$u$} at 322 310
\pinlabel {$v$} at 224 409
\endlabellist\plainweavecell} \qquad
{\labellist
\pinlabel {$X$} at 15 282
\pinlabel {$Y$} at 320 382
\endlabellist\plainweavekernel}\ec
\caption{Plain Weave} \label{plainweavefig}
\end{figure}

\begin{figure}[ht!]
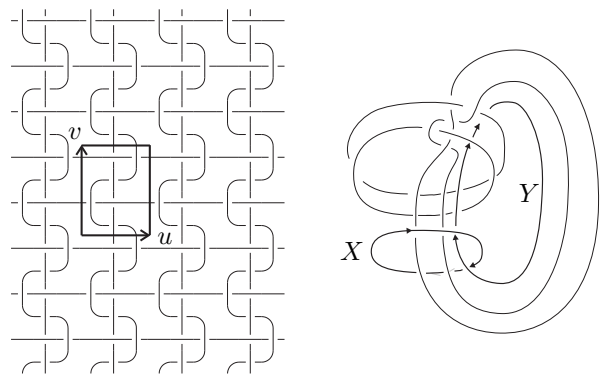

\bc
{\labellist
\pinlabel {$u$} at 333 359
\pinlabel {$v$} at 246 460
\endlabellist\lenoweavecell} \qquad
{\labellist
\pinlabel {$X$} at 15 282
\pinlabel {$Y$} at 320 382
\endlabellist\lenoweavekernel}\ec
\caption{Leno Weave} \label{lenoweavefig}
\end{figure}

\begin{figure}[ht!]
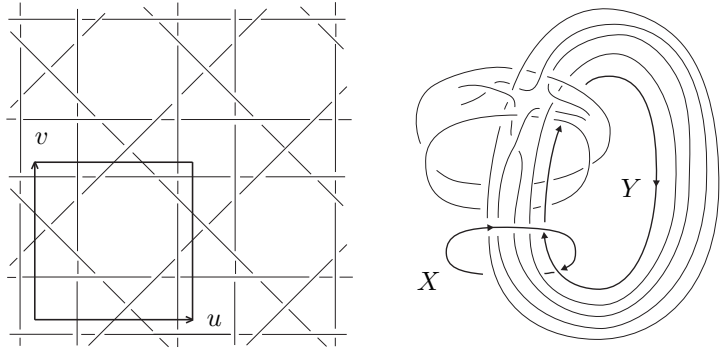

\bc
{\labellist
\pinlabel {$u$} at 261 34
\pinlabel {$v$} at 44 262
\endlabellist\fouraxialcell} \qquad
{\labellist
\pinlabel {$X$} at 16 210
\pinlabel {$Y$} at 212 300
\endlabellist\quatroaxialkernel}\ec
\caption{Multiaxial Weave} \label{quatroaxialfig}
\end{figure}

\begin{figure}[ht!]
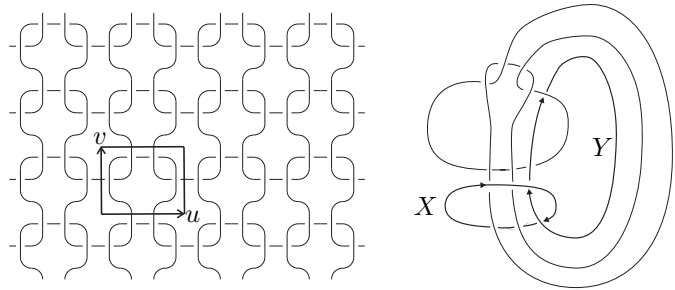

\bc
{\labellist
\pinlabel {$u$} at 323 316
\pinlabel {$v$} at 207 414
\endlabellist\singlejerseycell} \qquad
{\labellist
\pinlabel {$X$} at 15 282
\pinlabel {$Y$} at 320 382
\endlabellist\singlejerseykernel}\ec
\caption{Single Jersey} \label{singlejerseyfig}
\end{figure}

\newpage

\begin{figure}[ht!]
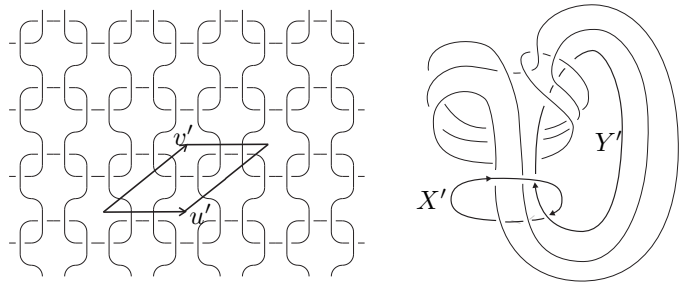

\bc
{\labellist
\pinlabel {$u'$} at 334 315
\pinlabel {$v'$} at 312 414
\endlabellist\singlejerseycellskew} \qquad
{\labellist
\pinlabel {$X'$} at 15 282
\pinlabel {$Y'$} at 320 382
\endlabellist\singlejerseykernelskew}\ec
\caption{Single Jersey with a different repeating cell} \label{singlejerseyskewfig}
\end{figure}

Different choices of repeating cell for a given fabric will give rise to different kernels. For
example, choosing the cell shown in Figure \ref{singlejerseyskewfig} for the single jersey
fabric gives the kernel shown.

\subsection{Fabric kernels}
 The original fabric can be recovered from any one of its
 kernels. Since the region between the two auxiliary components $X$ and $Y$ is topologically
 a thickened torus this means that the diagram on the torus can be recovered. The whole fabric
 can be reconstructed as a doubly periodic plane pattern by unwrapping the torus, - in effect
 constructing the inverse image of the diagram on the torus under the universal cover of the
 torus by the plane.

We adopt the name \emph{fabric kernel} for such links.
\begin{definition} A \emph{fabric kernel} is a link consisting of two distinguished unknotted
components $X$ and $Y$ which form a Hopf link and one or more further components representing
 the fabric strands.
\end{definition}
Any fabric kernel determines a fabric as above. Many fabric kernels will give rise to fabrics
which decompose into disconnected layers or strips, or can prove physically difficult to make.
The problem of identification of such structures which are not textiles from the traditional
point of view is a part of a wider problem of enumeration of all possible textile structures.
Some restricted subclasses corresponding to variants of traditional woven or knitted material
will be of particular interest to us, but much of our theoretical work will apply to very
general fabric links.

We show in the corollary to Theorem \ref{axial} how to use the classical multivariable
Alexander polynomial of a fabric kernel  to tell whether the resulting fabric contains any
closed components. Such \emph{chain mail} type of fabric is impossible to make using
traditional textile technology, and it is useful to be able to identify readily the
corresponding fabric kernels.

In section \ref{alexdata} we show how to compare the polynomials for different kernels for the
same fabric. The \emph{Alexander data}, which can be displayed as described in section
\ref{alexdata} on a diagram of the fabric, can be found from any given kernel, and allows the
polynomial of any other kernel to be readily calculated.

The Alexander data is then enough to allow comparison of different fabrics. While two different
fabrics can give rise to the same Alexander data the size of the kernels must be sufficiently
large for this to be possible, and an inventory of small fabric kernels using their polynomials
will give  a good means of listing the corresponding fabrics.  We expect that more detailed
comparisons of two fabrics will be best conducted by using kernels of approximately the same
size - measured for example by the number of crossings in a repeating cell.

The Alexander data  of a layered fabric has some characteristic features, described in section
\ref{layered}. Assuming that the layered fabric has no closed components then the multivariable
Alexander polynomial of any of its kernels must factorise non-trivially in a certain way,
described in Theorem \ref{layers}. This leads to a quick check for the maximum number of potential
layers in a fabric.

\section{The axial type of a strand in a fabric}
Suppose that we are given a fabric lying in a thickened plane, with  a discrete group $G$ of invariant translations generated by two independent elements $u$ and $v$. We shall assume that the normal direction $u\x v$  lies to the face side of the fabric.

Any strand of yarn in the fabric is invariant under some subgroup of $G$. This subgroup is either trivial, when the strand is a closed curve, or is infinite cyclic, generated by some $w\in G$. The element $w$ can be recognised, up to sign, as the smallest translation in $G$ that carries one point of the chosen strand to another point of the same strand. If the strand is oriented then we select $w$ so as to translate the strand in its preferred direction.

\begin{definition} The \emph{axial type} of a strand in a fabric is $0$ if the strand is closed and is otherwise a generator $w$ of the subgroup of translations which leaves the strand invariant.
\end{definition}

\begin{remark} Once we have chosen generators $u$ and $v$ for $G$ each axial type has the form $w=\alpha u+\beta v$ for some integers $(\alpha,\beta)$, determined up to an overeall sign if the strand is unoriented.
\end{remark}

For example, in the single jersey fabric shown in Figure \ref{singlejerseyfig} all fabric strands have axial type $u$. In fact all strands are equivalent under translations in $G$, and yield a single fabric component in the kernel.

In general there may be two or more inequivalent strands with the same axial type. For example in
Leno weave (Figure \ref{lenoweavefig}) there are two inequivalent strands with axial type $u$ and
two with axial type $v$, while the multiaxial fabric in Figure \ref{quatroaxialfig} has strands
of type $u,v,u+v$ and $u-v$.

\subsection{Axial type and kernels} \label{axialtype}

We can identify the axial type of the strands in a fabric by calculating certain linking numbers in a kernel for the fabric.

Assume that the kernel has been constructed using a choice of generators $u,v$ for the group $G$ of invariant translations. In the torus $T^2={\bf R}^2/G$ the lines in the directions of $u$ and $v$ become closed curves, $U$ and $V$ say, which generate the homology group $H_1(T^2)$. A strand in the fabric with axial type
$w=\alpha u+\beta v$ will become a closed curve $W$ in the thickened torus $T^2\x I$ which represents $\alpha U+\beta V$ in its homology group. Now the curves $X$ and $Y$ in the fabric kernel lie parallel to $U$ and $V$ respectively on either side of the standard thickened torus. They are oriented so that $lk(X,U)=0, lk(X,V)=1, lk(Y,U)=1, lk(Y,U)=0$. Since $W=\alpha U+\beta V$  in the homology of $T^2\x I\cong S^3-(X\cup Y)$ we get $lk(X,W)=\alpha lk(X,U)+\beta lk(X,V)=\beta$ and similarly $lk(Y,W)=\alpha$.

Consequently in a fabric kernel $L=X\cup Y\cup T_1\cup\ldots\cup T_k$ each oriented component $T_i$ represents a family of translation-equivalent strands in the fabric of axial type $b_i u +a_i v$, where $a_i=lk(X,T_i)$ and $b_i=lk(Y,T_i)$.  In particular the strands in the fabric corresponding to $T_i$ are \emph{closed curves} if and only if $a_i=b_i=0$.

\section{The multivariable Alexander polynomial}

For a link $L\subset S^3$ with $r>1$ oriented components $X_1\cup\ldots\cup X_r$ the group
$H_1(S^3-L)$ is free abelian of rank $r$. It has distinguished generators $x_1,\ldots,x_r$
represented by oriented meridians of the components. The multivariable Alexander polynomial
$\Delta_L$ is an element of the integer group ring ${\bf Z}[H_1(S^3-L)]$ and is thus a Laurent
polynomial in ${\bf Z}[x_1^{\pm1},\ldots,x_r^{\pm1}]$. It is defined up to a unit in this ring,
and thus up to multiplication by a signed monomial $\pm x_1^{a_1}\cdots x_r^{a_r}$.

With some care an absolute version of the polynomial can be defined, as in \cite{Murakami}. The
Torres symmetry condition (\ref{TF}) shows that $$\Delta_L(x_1^{-1},\ldots,x_r^{-1})=(-1)^r
x_1^{m_1}\ldots x_r^{m_r}\Delta_L(x_1,\ldots,x_r),$$ for some integers $m_1,\ldots,m_r$.
Multiplying $\Delta_L$ by the monomial $x_1^{\ha m_1}\ldots x_r^{\ha m_r}$ gives Murakami's
absolute version, up to sign, with the property that
$$\Delta_L(x_1^{-1},\ldots,x_r^{-1})=(-1)^r
\Delta_L(x_1,\ldots,x_r).$$ In this format we need to use half-integer powers $x_i^\ha$ of some of the
variables $\{x_i\}$.

We make use of the absolute version in presenting the Alexander data, but for many of the properties and calculations it is enough to allow a monomial multiple, so as to avoid negative powers of the variables.

In the case of a single component link, in other words a knot $K$, we write
$\overline\Delta_K\in {\bf Z}[t^{\pm1}]$ for its classical Alexander polynomial, and in this
paper we use the non-standard notation $\Delta_K$ to denote the rational function
$\Delta_K=\frac{1}{1-t}\overline\Delta_K$. In this way we can give formulae for the behaviour
of the Alexander polynomials of related knots and links in a uniform way, without having to
treat the single component case,  $r=1$,  separately.  The formulae can be most consistently
handled in the context of Reidemeister torsion, where an excellent recent account is given in
\cite{Turaev}. Here we give a summary of the properties needed, following the constructions of
Fox and Torres based on Fox's free differential calculus \cite{TF, FoxV}.

\subsection{Deletion of components}
For an oriented link $L$ with components $X_1\cup\ldots\cup X_r$ the curve $X_i$ represents the
monomial $\ds \prod_{j\ne i}x_j^{l_{ij}}$ as an element of $H_1(S^3-(L-X_i))$, where $l_{ij}= lk(X_i,X_j)\in
{\bf Z}$ is the linking number  of   $X_i$ and $X_j$. We write
\be<X_i>=\prod_{j\ne i}x_j^{l_{ij}} \label{curvemonomial}\ee for this monomial. A special case of the Torres-Fox
satellite formula relates the invariants $\Delta_L$ and $\Delta_{L-X_i}$.
\begin{theorem}[Torres] \label{TF} $$\Delta_L|_{x_i=1}=(1-<X_i>)\Delta_{L-X_i}.$$
\end{theorem}

\begin{remark} If the component $X_i$ has a non-zero linking number with \emph{at least one} of
 the other components then $<X_i>\ne1$. It is thus possible to recover the Alexander polynomial
 of the sublink $L-X_i$ where $X_i$ is removed, starting from the polynomial $\Delta_L$ of the
 whole link, unless $X_i$ has linking number $0$ with \emph{all} the other components of $L$.

For the absolute versions of the invariants in Theorem \ref{TF}  the factor \break ${<X_i>^\ha
-<X_i>^{-\ha}}$ appears in place of $1-<X_i>$.
\end{remark}

\subsection{Closed components in a fabric}
We start from a fabric kernel $L$. Then $L=X\cup Y\cup T_1\cup\ldots\cup T_k$, where $X$ and
$Y$ are the distinguished face and back curves forming a Hopf link. The complement of  $X \cup
Y$ forms a standard thickened torus containing the remaining components  $T_1,\ldots, T_k$,
termed the \emph{fabric components} of $L$. These curves correspond to the strands of yarn in
the fabric resulting from unwrapping $L$.

The curve $T_i$ in the fabric kernel unwraps to give closed components in the covering fabric
if and only if it has linking number $0$ with both of the auxiliary curves $X$ and $Y$.

\begin{theorem}\label{axial}
Let $L=X\cup Y\cup T_1\cup\ldots\cup T_k$ be a fabric kernel.
Write $a_i=lk(T_i,X)$ and $b_i=lk(T_i,Y)$. Then its Alexander polynomial $\Delta_L(x,y,t_1,\ldots,t_k)$ satisfies $$\Delta_L(x,y,1,\ldots,1)=\prod_{i=1}^k (1-x^{a_i}y^{b_i})\in {\bf Z}[x^{\pm1},y^{\pm1}].$$
\end{theorem}
\begin{proof} After setting all $t_j=1$ we have $<T_i>=x^{a_i}y^{b_i}$ by equation (\ref{curvemonomial}). Repeated use of
theorem \ref{TF}, suppressing the components $T_1, \ldots , T_k$ in turn,  gives $$\Delta_L(x,y,1,\ldots,1)=\prod_{i=1}^k (1-x^{a_i}y^{b_i})\Delta_{X\cup Y}.$$ Since the remaining link $X\cup Y$ is the Hopf link, whose Alexander polynomial is $1$, the result follows.
\end{proof}

\begin{corollary}
There is a closed component in the covering fabric of $L$ if and only if $\Delta_L(x,y,1,\ldots,1)=0.$
\end{corollary}
\begin{proof}
The Laurent polynomial
$\prod_{i=1}^k (1-x^{a_i}y^{b_i})$
is equal to $0$  in ${\bf Z}[x^{\pm1},y^{\pm1}]$ if and only if $a_i=b_i=0$ for some $i$.
\end{proof}

\begin{remark} The element $<T_i>=x^{a_i}y^{b_i}$ represents the homology class of $T_i$ in the thickened torus $S^3-(X\cup Y)$. Written additively this is $b_iU+a_iV$, where $U$ and $V$ are determined by the generators $u$ and $v$ as in subsection \ref{axialtype}.
The axial type $b_i u+a_i v$ of the corresponding strands in the fabric, and indeed the number of inequivalent strands of each axial type,  can then be read off immediately from $\Delta_L$, so long as there are no closed strands in the fabric.

To determine the axial types using Theorem \ref{axial} it is enough to be given the polynomial $\Delta_L$ as a function of $x,y$ and $t$
only, where all the yarn variables $t_i$ have been set equal to $t$.
 If there are no closed components in the fabric we can then recover
the number $k$ of  yarn components in the kernel, and the axial types of the corresponding fabric strands,  from the factors in the evaluation with $t=1$.
\end{remark}

 In a
traditional woven fabric there are just two axial types, corresponding to the warp and weft
directions, while   a traditional knitted fabric has only one type.

\begin{figure}[ht!]
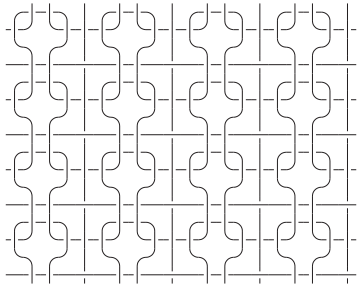

\bc
\singlejerseywarpweft  \ec
\caption{Single Jersey with warp and weft inlays} \label{singlejerseywarpweftfig}
\end{figure}

More sophisticated woven fabrics with multiaxial types have recently been introduced, \cite{Dow70}, such as
the fabrics shown in Figure  \ref{quatroaxialfig} and \ref{triaxialweavefig}. Knitted fabrics too may have warp and
weft type inserts, as in Figure \ref{singlejerseywarpweftfig},  leading to a multiaxial fabric.

From this basic analysis of the polynomial of a kernel we can identify immediately the axial
types, although a detailed discussion of whether a fabric is woven or knitted
in the appropriate sense may not be available from the polynomial alone.  Of course if the
kernel has more than $2$ different axial types then we can certainly conclude that the
fabric is not a traditional woven structure, while equally a traditional knitted structure can
be excluded unless there is a single axial type.

We give here the polynomials for the kernels of several fabrics.

\begin{example}
\noindent Single jersey with a closed thread around the fabric thread

$$\Delta_L=(tx-t-x)(t+x-1)(t-1)(y-1).$$

\end{example}

\begin{example}

\noindent Single jersey with a trefoil added to the fabric thread

$$\Delta_L=(tx-t-x)(t+x-1)(t^2-t+1)(y-1).$$

\end{example}
\begin{example}
\noindent Single jersey with warp and weft inlays, as in Figure \ref{singlejerseywarpweftfig}:
$$\Delta_L=(t+x-1)(tx-t-x)(px-1)(ey-1)(y-1),$$ where the jersey strand has variable $t$, the warp strand has variable $e$ and the weft has variable $p$.
\end{example}

\begin{example}

By way of contrast, the Alexander polynomial for the kernel of the chain mail pattern shown in
Figure \ref{chainmailfig} is \be\label{chainmailpoly}\Delta_L= (x-y)(1-xy)(1-t).\ee

\begin{figure}[ht!]
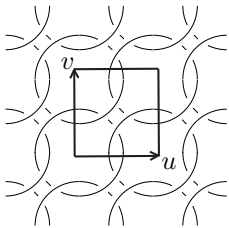

\bc
{\labellist
\pinlabel {$u$} at 221 180
\pinlabel {$v$} at 123 277
\endlabellist\chainmail} \ec
\caption{Chain Mail} \label{chainmailfig}
\end{figure}

This is one of the fabric kernels which does not correspond to a physically constructible textile fabric, on account of the closed components. The substitution $t=1$ for the fabric variables leads to the value $0$, as required by Theorem \ref{axial}.
\end{example}

\subsection{The Alexander data for a fabric} \label{alexdata}
We now present a normalised form of the Alexander polynomial of a kernel giving rise to an array of polynomials in the fabric variables which are independent of the choice of generators $u$ and $v$. The result is a set of invariants for the fabric which readily determine the Alexander polynomial for \emph{any} choice of kernel.

Start from the Alexander polynomial $\Delta_L(x,y,t_1,\ldots,t_k)$ of a kernel $L$ determined by a choice of generators $u$ and $v$ of the group of translations of the fabric.

Write $t_i=s_i^2$ and set \be U=y\prod s_i^{a_i}, V=x\prod s_i^{b_i} \label{alexandervariable}\ee
where $a_i=lk(T_i,X)$ and $b_i=lk(T_i,Y)$. The integers $a_i$ and $b_i$ can be easily read off
from the unit cell determined by $u$ and $v$ using the fact that $T_i$ crosses  the edge $u$
algebraically $a_i$ times in the direction of $v$ and  crosses the edge $v$ algebraically $b_i$
times in the direction of $u$.

The \emph{Alexander data} of the fabric consists of the coefficients of $\Delta_L$ when written as a polynomial in $U$ and $V$. When written in this way we say that the Alexander polynomial $\Delta_L$ is in \emph{data form}.

The Alexander data  coefficients are Laurent polynomials in the fabric variables $\{s_i\}$.   In Theorem \ref{alexanderdata} we show how they depend essentially on the fabric only and not on the choice of kernel.

Write $$\Delta_L=\sum L_{\alpha,\beta}U^\alpha V^\beta$$ and set $T_w=L_{\alpha,\beta}$ for each $w=\alpha u+\beta v\in G$. In this way we associate a Laurent polynomial $T_w(s_1,\ldots,s_k)$ to each element $w$ in the translation group $G$.

Since the Alexander polynomial $\Delta_L$ is often defined only up to multiplication by a signed monomial in $x,y,t_1,\ldots,t_k$ there is a measure of ambiguity in the definition of $T_w$. Apart from this ambiguity, which can be eliminated up to an overall sign by use of the Torres symmetry conditions, the elements $T_w$ are independent on the choice of kernel, and are termed collectively the \emph{Alexander data}  of the fabric. The independence is made precise in the following invariance theorem.

\begin{theorem}
 \label{alexanderdata}
The Alexander data of a fabric is independent of the choice of generators for $G$, up to
an overall multiplication by a monomial in $s_1,\ldots,s_k$ and translation by an element
of $G$. When we use generators $u'$ and $v'$ in place of $u$ and $v$ to form a kernel $L'$, then the resulting polynomials $T'_w, w\in G$ satisfy the equation
$$T'_w=\pm m T_{w+g},$$ for some $g\in G$ and monomial $m$ in $s_1,\ldots,s_k$ independent of $w$.
\end{theorem}

\begin{proof}
  Our choice $u,v$ of generators for the group of translations $G$ determines the repeating
  cell used to construct the kernel.  Lines in the plane in the directions of $u$ and $v$
  become closed curves $U$ and $V$ in the quotient torus $T^2 ={\bf R}^2/G$. The face torus
  is embedded as the boundary of a neighbourhood of the face curve $X$ and the back torus as
  the boundary of a neighbourhood of the back curve $Y$, where $X$ and $Y$ are parallel to
  the closed curves $U$ and $V$ respectively. In this embedding the face torus contains
  curves $m_X$ and $l_X$, which are respectively the meridian and longitude of the curve $X$.
  The back torus correspondingly contains $m_Y$ and $l_Y$, giving the longitude and meridian
  curves as $l_X=U\x \{1\}$, $l_Y=V\x\{0\}$, while $m_X=V\x\{1\}$ and $m_Y=U\x\{0\}$.

In terms of the generators $x,y,t_1,\ldots, t_k$ for the homology $H_1(S^3-L)$ we have \beqn U\x\{0\}=y,  &&U\x\{1\}=<X>=y\prod t_i^{a_i},\\
V\x \{1\}=x, && V\x\{0\}=<Y>=x\prod t_i^{b_i},\eeqn where $a_i=lk(T_i,X)$ and $b_i=lk(T_i,Y)$.  The   variables $U=y\prod s_i^{a_i}$ and $V=x\prod s_i^{b_i}$ introduced in equation (\ref{alexandervariable})  thus represent  an `average' of the homology of the curves $U$ and $V$ on the front and back faces of the torus. When the homology $H_1(S^3-L)$ is written additively we have \beqn 2U&=&U\x\{0\}+U\x\{1\}\\ 2V&=&V\x\{0\}+V\x\{1\}.\eeqn

A different choice of generators $u',v'$ for $G$ gives rise to a different repeating cell,
and a different kernel for the fabric.  Now the complement of a kernel in $S^3$ is the
complement of the fabric quotient in the thickened torus $T^2\x I$. This complement is unchanged by the
 new choice of generators for $G$, so that $S^3-L'\cong S^3-L$.  What \emph{does} change is the curves on the
face and back tori that correspond to the longitudes and meridians of the face and back curves in the two kernels. We then know that $\Delta_L'=\Delta_L$, when the variables are changed appropriately.  We need to compare the elements $x,y,t_1,\ldots,t_k$ and $x',y',t_1,\ldots,t_k$ as elements of $H_1(S^3-L)=H_1(S^3-L')$. The comparison is very straightforward when we use the data form of the two polynomials.

\begin{lemma}\label{datavariablechange}
We can suppose that
$$ u'=pu+qv, v'=ru+sv$$ with $p,q,r,s\in{\bf Z}$ and $ps-qr=1$. Then $$ U'=pU+qV, V'=rU+sV, $$ when written additively as elements of $H_1(S^3-L)$.
\end{lemma}

\begin{proof}
The curves $U', V'$ on $T^2={\bf R}^2/G$ corresponding to the vectors $u',v'$ satisfy the equations $$ U'=pU+qV,\quad V'=rU+sV, $$ as elements of $H_1(T^2)$.

Hence on the boundary of the thickened torus we have $$ U'\x\{i\}=pU\x\{i\}+qV\x\{i\}, \quad V'\x\{i\}=rU\x\{i\}+sV\x\{i\}, $$ for $i=0,1$.

Adding these equations for $i=0$ and $i=1$ in $H_1(S^3-L)$ we then have $$ 2U'=2(pU+qV), 2V'=2(rU+sV). $$

\end{proof}

The multivariable Alexander polynomial $\Delta_L$  of a fabric kernel
$L=X\cup Y\cup T_1\cup\ldots\cup T_k$ is an element of the integer group ring
${\bf Z}[H_1(S^3-L)]$, and hence is a Laurent polynomial in variables $x,y, t_1,\ldots,t_k$,
since the homology group can be generated by meridians of the components. The element
$\Delta_L$ depends only on the fundamental group of $S^3-L$, which is also the fundamental
group of the complement of the fabric quotient in the thickened torus $T^2\x I$. Hence
if $L'$ is another kernel for the same fabric it will have the same Alexander polynomial,
since $\pi_1(S^3-L')=\pi_1(S^3-L)$, written now in terms of different generators of the
group $H_1(S^3-L)$.

In terms of the Alexander polynomial we can relate $\Delta_{L'}$ to $\Delta_L$, up to multiplication by a signed monomial, using the multiplicative version $U'=U^pV^q$, $V'=U^rV^s$ of Lemma \ref{datavariablechange}.

The Alexander data polynomials $T'_w$ defined using the kernel $L'$ are given by the equation $\Delta_{L'}=\sum T'_w (U')^{\alpha'}(V')^{\beta'}$, where $w=\alpha'u'+\beta'v'$.  Making the substitution for $U'$ and $V'$ gives $$\Delta_L=\Delta_{L'}=\sum T'_w U^{(\alpha' p+\beta' r)}V^{(\alpha' q+\beta' s)}=\sum T'_w U^\alpha V^\beta,$$ with $\alpha=\alpha' p+\beta' r, \beta =\alpha' q+\beta' s$. Hence, in terms of the Alexander data defined using $L$ we have $T'_w=T_g$ with $g=\alpha u+\beta v$.  Now $g=(\alpha' p+\beta' r)u+(\alpha' q+\beta' s)v=\alpha'u'+\beta'v'=w$, giving the invariance result that $T'_w=T_w$, up to the overall ambiguity of the definition of the Alexander polynomial.
\end{proof}

\subsection{Changing the choice of the repeating cell}
 Given the Alexander data of a fabric, calculated using any one choice of generators for the translation group, Theorem \ref{alexanderdata} shows how  the multivariable
 Alexander polynomial for any other kernel of the same fabric can be found.

We can now compare the Alexander data from two fabrics quite readily. We just need one
choice of unit cell for each  to provide the Alexander data. Having put the two sets of
data onto a plane we can see that if the fabrics themselves are affinely equivalent, in
other words related by a translation and a linear transformation, then one set of data
will transform to the other  by a similar affine transformation $\varphi$ carrying the
coefficient polynomials $\{T_w\}$ for one fabric to the polynomials $\{T_{\varphi(w)}\}$ of
the other.  Fabrics with inequivalent  Alexander data are then inequivalent.

In particular the number and affine location of the non-zero polynomials in the Alexander
data is a simple invariant of the fabric. This could be refined, without presenting the whole data, to include the degrees of each fabric variable at the plane locations.

\subsection{Presenting the Alexander data}

Superimpose a grid on the fabric by choosing a point $O$ as origin in the plane, and label
the translates of $O$ by elements of the group $G$. At the translate of $O$ by $w\in G$ place the
polynomial $T_w$. This display of the polynomials $T_w$ is uniquely determined by the choice
of origin, where the polynomials themselves have been normalised, up to an overall sign, by
multiplying all polynomials by a monomial in the variables $\{s_i=t_i^\ha\}$ so as to respect
the Torres symmetry conditions.

Thus to find the Alexander data polynomials $T_w$ we choose any kernel $L$ and put $\Delta_L$ into data form by setting $x=V/\prod s_i^{b_i}$ and $y=U/\prod s_i^{a_i}$ as in equation (\ref{alexandervariable}).  Then $T_{\alpha u+\beta v}$ is the coefficient of $U^\alpha V^\beta $ in the resulting polynomial.

Conversely, given $u$ and $v$ we can recover the data form of $\Delta_L$ in terms of $U$ and $V$ immediately from the Alexander data. Equation (\ref{alexandervariable}), using the numbers $a_i$ and $b_i$ calculated from the choice of $u$ and $v$, then gives $\Delta_L$  in terms of $x$ and $y$.

\begin{example}\label{jerseyexample}

We may compare the single jersey data derived from  the different
choices of unit cell shown in Figures \ref{singlejerseyfig} and \ref{singlejerseyskewfig}.

From a braid presentation of the single jersey kernel $L$ shown in Figure \ref{singlejerseyfig}
we can calculate its multivariable Alexander polynomial
$$\Delta_L=(1-y)(1-x-t)(x+t-tx),$$ where the variable $t=s^2$ is the single fabric variable.
Setting $t=1$ gives $y-1$ up to monomial multiples, confirming that  the fabric has a single
axial direction $u$. From the intersections with the cell edges we have the Alexander data substitution $U=y$ and $V=xs$.

A similar braid based calculation for the polynomial $\Delta_{L'}$ of the kernel in Figure \ref{singlejerseyskewfig} gives
$$\Delta_{L'}=(1-y')(x'-tx'+ty')(x'-y'+ty')$$ up to a signed monomial multiple.  Again the intersections with the cell edges show that the data substitution is $U'=y', V'=x's$. Now this
kernel is related to $L$ by the choice of generators $u'=u,v'=u+v$. Consequently the substitution $U'=U, V'=UV$ will convert the data form for $L'$ into the data form for $L$, up to a monomial multiple.

We can exhibit the Alexander data  starting from either data form. On the one hand we have
$$\Delta_L=(1-U)(1-s^2-V/s)(s^2+(s^{-1}-s)V),$$ while
$$\Delta_{L'}=(1-U')((s^{-1}-s)V'+s^2U')(V'/s-(1-s^2)U')=U^2\Delta_L$$ after the substitution.

 Multiplication  by $s^{-1}$ gives the Torres symmetric form of the data displayed in matrix form below.

\begin{figure}[ht!]
\bc
$\begin{array}{cc}
s^{-3}-s^{-1}&s^{-1}-s^{-3}\\[4mm]
-s^2+3-s^{-2}&s^2-3+s^{-2}\\[4mm]
s^3-s&s-s^3\end{array}$
\ec \caption{Single Jersey (face) data} \label{singlejerseydatafig}
\end{figure}

\end{example}

Here and in the examples below we present the Alexander data as a matrix of polynomials based on the stated choice of $u$ and $v$. The rows of the matrix correspond to translations by multiples of $u$, read from left to right,  and the columns to multiples of $v$, read from the bottom up.

\begin{example}\label{jerseybackexample}

In the previous example we have used the single jersey fabric viewed from its traditional face side. When viewed from the back, as in Figure \ref{singlejerseybackfig}, the polynomial for the kernel is $\Delta_L=(1-y)(t^2x-tx+1)(t^2x-t+1)$.
\begin{figure}[ht!]
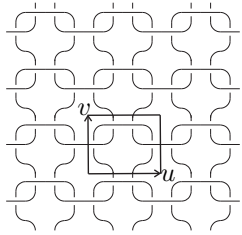

\bc
{\labellist
\pinlabel {$u$} at 321 360
\pinlabel {$v$} at 258 410
\endlabellist\singlejerseyback} \ec
\caption{Single Jersey, from the back} \label{singlejerseybackfig}
\end{figure}
 This gives the normalised Alexander data shown in Figure \ref{singlejerseybackdatafig}.
\begin{figure}[ht!]
\bc
$\begin{array}{cc}
s^3-s&s-s^3\\[4mm]
-s^2+3-s^{-2}&s^2-3+s^{-2}\\[4mm]
s^{-3}-s^{-1}&s^{-1}-s^{-3}\end{array}$
\ec
\caption{Single Jersey (back) data} \label{singlejerseybackdatafig}
\end{figure}

\end{example}

\begin{example}\label{plainweaveexample}
The Alexander polynomial for the kernel $L$ of plain weave  corresponding to the choices of $u$ and $v$ in Figure \ref{plainweavefig} is
\beqn \Delta_L&=& -1+2px+2ey-p^2 x^2-e^2y^2 +((1-e)^2(1-p)^2 -4ep)xy \\ && +2p^2ex^2y +2pe^2xy^2-e^2p^2x^2y^2,\eeqn where the two fabric strands in the direction of $u$ have the meridian variable $p$ and those in the direction $v$ have variable $e$. Write $e=a^2$ and $p=b^2$.   Then setting $U=y/a^2=y/e$ and $V=x/b^2=x/p$, as required by equation (\ref{alexanderdata}), gives the data form  $$\Delta_L= -1+2V+2U-V^2-U^2 +((1-e)^2(1-p)^2-4ep) UV/ep +2UV^2 +2U^2V-U^2V^2.$$ This leads to the data shown in Figure \ref{plainweavedata}.

\begin{figure}[ht!]
\bc
$ \begin{array}{ccc}
-1&2&-1\\[4mm]
2&(a-a^{-1})^2(b-b^{-1})^2-4&2\\[4mm]
-1&2&-1
\end{array}$
\ec \caption{Plain Weave data} \label{plainweavedata}
\end{figure}
\end{example}

\newpage
\begin{example}\label{lenoweaveexample}

Leno weave, which was presented in Figure \ref{lenoweavefig}, has four threads in the minimal
unit cell; two non-equivalent threads are running in the $u$ direction and two in the $v$ direction.

The Alexander data for Leno weave is presented in Figure \ref{lenoweavedata}, where \beqn P&=&(a^{-2}-1+a^2)(b^{-1}-b)^2+2\\ Q&=&2(a^{-1}-a)^2(b^{-1}-b)^2+2(a^{-2}-a^2)^2-(a^{-2}b^{-1}+a^2b)^2+2.\eeqn

\begin{figure}[ht!]
\bc
$ \begin{array}{ccc}
-1&2&-1\\[4mm]
P&Q&P\\[4mm]
-1&2&-1
\end{array}$
\ec \caption{Leno weave data} \label{lenoweavedata}
\end{figure}

\end{example}

\begin{example}\label{chainmailexample}
The Alexander data for chain mail is   shown in Figure \ref{chainmaildata}.
It is based on the choice of generators in Figure \ref{chainmailfig}, and the corresponding Alexander polynomial in equation (\ref{chainmailpoly}). In this case the data form is given simply by $U=y, V=x$.
\begin{figure}[ht!]
\bc
$ \begin{array}{ccc}
0&s-s^{-1}&0\\[4mm]
s^{-1}-s&0&s^{-1}-s\\[4mm]
0&s-s^{-1}&0
\end{array}$
\ec \caption{Chain Mail data} \label{chainmaildata}
\end{figure}
\end{example}
\newpage

\begin{example}\label{triaxialexample}
The triaxial weave shown in Figure \ref{triaxialweavefig} has a lattice of symmetries with a $3$-fold symmetry when the strings are oriented in the three axial directions, $u$, $v$ and $w=-u-v$  as shown.

\begin{figure}[ht!]
\bc
{\labellist
\pinlabel {$u$} at 326 399
\pinlabel {$v$} at 253 537
\pinlabel {$w$} at 170 405
\endlabellist\triaxial} \ec
\caption{Triaxial weave} \label{triaxialweavefig}
\end{figure}

The Alexander data, when presented on this lattice,   exhibits the same symmetry very nicely. This is shown in Figure \ref{triaxialweavedata}, where the meridian variables of the strings are $a^2,b^2$ and $c^2$ respectively, and $P=(a-a^{-1})(b-b^{-1})(c-c^{-1})$.

\begin{figure}[ht!]
\bc
{\labellist
\pinlabel {$a$} at 103 531
\pinlabel {$b$} at 291 420
\pinlabel {$c$} at 291 640

\pinlabel {$-a^{-1}$} at 350 531
\pinlabel {$-b^{-1}$} at 164 640
\pinlabel {$-c^{-1}$} at 164 420
\pinlabel {$P$} at 223 531

\endlabellist\triaxialdata} \ec
\caption{Triaxial weave} \label{triaxialweavedata}
\end{figure}
\end{example}

\subsection{Multiple repeating cells}
For comparison purposes it can be useful to look not just at a minimal area repeat of a
pattern in a fabric. This means that we use a repeating cell whose edges are determined
by two vectors $u',v'$ in the group $G$ of translations which generate a proper subgroup
of $G$. It is always possible to choose generators $u,v$ for $G$ and $u',v'$ for a subgroup $H$
of $G$ so that $u'$ and $v'$ are multiples of $u$ and $v$ respectively. We can insist further
that $u'=su$ and $v'=rv$ where $r$ is a multiple of $s$. Then every vector $w'$  in $H$ is a
multiple $w'=sw$ of a vector $w$ in $G$, and in many cases we will have $s=1$.

The relation between the Alexander data for the fabric using the full group $G$, corresponding
to a minimal choice of repeating cell, and the data derived from using the translations in $H$
only (and thus only considering a larger repeating area) can be determined by comparing the
Alexander polynomials of the kernels $L$ and $L^{(r)}$ which come from the choice of
translations $u,v$ and $u,w=rv$ respectively.

These polynomials are related by a more general result of Salkeld  \cite{Salkeld} connecting
the multivariable polynomials for a link $L$ with a distinguished unknotted component $X$ and
the link $L^{(r)}$ given by taking the $r$-fold cyclic cover of $L$ branched over $X$. The two
kernels above are related in this way, so Salkeld's result applies.

There is an $r$-fold covering map $p$ from the complement of $L^{(r)}$ to the complement
of $L$, inducing a map $p^*$ on the first homology group. Salkeld's theorem expresses
$p^*(\Delta_{L^{(r)}}) $ in terms of $\Delta_L$ as follows.

\begin{theorem}[Salkeld] Write $x$ for the meridian variable corresponding to the component
$X$ in $L$. Then
$$p^*(\Delta_{L^{(r)}}) =\prod \Delta_L(\zeta x),$$ as $\zeta$ runs through all $r$th roots
of $1$.
\end{theorem}

\begin{corollary}

The data form of $\Delta_{L^{(r)}}$, in terms of $U$ and $W$, corresponding to the translations $u$ and $w=rv$, satisfies the equation $$p^*(\Delta_{L^{(r)}}) =\prod \Delta_L(\zeta V),$$ as $\zeta$ runs through all $r$th roots
of $1$, where $W=V^r$.
\end{corollary}
\begin{proof}
In the data form for  $ \Delta_L$ we have $V=x\prod s_i^{a_i}$, while in $\Delta_{L^{(r)}}$ we have $W=X(\prod s_i^{a_i})^r$, and $p^*(X)=x^r$. Now $ \Delta_L(\zeta x)$,  given by replacing $x$ with $\zeta x$, can be found from the data form by replacing $V$ with $\zeta V$.
\end{proof}

Salkeld's theorem does not give the complete polynomial for the covering link, in the case when
there is more than one component of $L^{(r)}$ covering some component of $L$, since
$p^*(\Delta_{L^{(r)}})$ has the same variable for all components which project to the same
component in $L$. In the case of fabric kernels we would expect to make this restriction in any
case, as the different components of  $L^{(r)}$ with the same image must be equivalent strands
in the fabric under the full group $G$ of translations. We would only use a different variable
if they were to be considered different in the fabric, and in this case a translation carrying
one to the other would not be part of the group $G$. The meridian for the branch curve in
$L^{(r)}$ projects to $x^r$ in $L$, and so the polynomial $\Delta_{L^{(r)}}$, with equivalent
components in the fabric having the same labels, can be recovered from $\Delta_L$.
 Conversely, given $\Delta_{L^{(r)}}$ we can recover $\Delta_L$ as one factor, corresponding
 to the root $\zeta=1$, when the branch curve meridian is replaced by $x^r$.

If $\Delta_L$ factorises, then $\Delta_{L^{(r)}}$ can be found as a product using  the same
operation on each factor. Indeed the operation of passing from
$\Delta_L$ to $p^*(\Delta_{L^{(r)}})$ can be regarded formally as replacing the roots
(for $x$) of $\Delta_L$ by their $r$-th powers;   the coefficients of powers of $x$
in $p^*(\Delta_{L^{(r)}})$ are integer polynomials in the coefficients of $\Delta_L$.

\begin{example} \label{jerseydoublingexample}
We can find the data form of the Alexander polynomial for single jersey doubled horizontally,
using the cell  with sides $w=2u$ and $v$ in terms of the generators $u$ and $v$ shown in Figure
\ref{singlejerseyfig} for the translation group. This is given from Salkeld's theorem by taking
the 2-fold cover of the kernel branched over the curve $Y$. The result, as a polynomial in $W$
and $V$,  will be the polynomial $P(W,V,t)=\Delta_L(U,V,t)\Delta_L(-U,V,t) $ with $W=U^2$. Using
the data form $\Delta_L=(1-U)(1-s^2-V/s)(s^2+(s^{-1}-s)V)$  we get
\beqn P&=&(1-U)(1+U)(1-s^2-V/s)^2(s^2+(s^{-1}-s)V)^2\\
&=&(1-W)(1-s^2-V/s)^2(s^2+(s^{-1}-s)V)^2.\eeqn

The polynomial $Q$ for single jersey doubled vertically, corresponding to $u$ and $w=2v$ has $Q(U,W,t)=\Delta_L(U,V,t)\Delta_L(U,-V,t) $ with $W=V^2$, giving the data form
$$Q=(1-U)^2((1-t)^2-W/t)(t^2-(1-t)^2 W/t),$$ where both fabric strands have the variable $t$. This converts to the Alexander polynomial by setting $U=y, W=xs^2=xt$, to get $$Q=(1-y)^2((1-t)^2-x)(t^2-(1-t)^2 x).$$ A more refined version of the Alexander polynomial, using different variables for the two strands, can be calculated directly from a kernel as $$(1-y)^2((1-t_1)(1-t_2)-x)(t_1t_2-(1-t_1)(1-t_2)x).$$

\end{example}

\begin{example}\label{ribexample}
The Alexander polynomial for $1\x1$ rib, which alternates face and back loops
 in the  $u$
direction, is
$$\Delta_L=(tx-t-x)(t+x-1)(t^2x-tx+1)(t^2x-t+1)(y-1).$$
This polynomial is a mix of the polynomials for the  face and back kernels of single jersey, and can be compared with the  multiple cell case of single jersey in the previous example, where two face loops are repeated in the $u$ direction. This was calculated above to be $$\Delta_L=(tx-t-x)^2(t+x-1)^2(y-1).$$ A repeat of the back loop in the $u$ direction gives $$\Delta_L=(t^2x-tx+1)^2(t^2x-t+1)^2(y-1).$$ 
These two polynomials are inequivalent to the polynomial for $1\x 1$ rib.

\noindent For the general case of $m \times n$ rib, which is a combination of $m \geq 0$ face and
$n \geq 0$ back loops in the $u$ direction, we have

$$\Delta_L=((tx-t-x)(t+x-1))^m((t^2x-tx+1)(t^2x-t+1))^n(y-1).$$

Indeed this polynomial is the same irrespective of the order in which the face and back loops
occur, so that the standard $2\x2$ rib gives the same data as the $1\x1$ rib repeated twice in
the $u$ direction. We have tried comparing the polynomials in which alternate rows have variables
$t_1$ and $t_2$, to see if this is enough to distinguish these fabrics, but again the results are
the same.

\end{example}

\begin{example}\label{garterexample}

\noindent Alternate knit and purl rows, sometimes known as \emph{garter stitch}.
This is a combination of one face and one back loop in the $v$ direction, and has Alexander polynomial
$$\Delta_L=(1+(1-t_1)(1-t_2)x)((1-t_1)(1-t_2)+t_1 t_2x)(y-1)^2.$$

It can be compared with the polynomial for plain knit or purl fabric, calculated  with different variables $t_1$ and $t_2$ for alternate rows.

\noindent For two face loops in the $v$ direction the polynomial is

$$\Delta_L=(t_1 t_2-(1-t_1)(1-t_2)x)((1-t_1)(1-t_2)-x)(y-1)^2.$$

\noindent For two back loops in the $v$ direction the polynomial is

$$\Delta_L=({t_1^2}{t_2^2}x-t_1 t_2+t_1+t_2-1)({t_1^2}{t_2^2}x-{t_1^2}t_2x-{t_2^2}t_1x-t_1 t_2x-1)(y-1)^2.$$

Replacing $x$ by $1/t_1 t_2x$ converts the purl to the plain version here, and leaves the  polynomial of the mixed fabric unchanged. It is possible to use the different signs in the brackets to distinguish the mixed fabric with alternately knit and purl rows from the plain fabric viewed from either side, by setting $t_1=t_2=-1$.
\end{example}

\begin{example} We can use Salkeld's theorem as in Example \ref{jerseydoublingexample} to calculate the polynomials for plain weave when doubled in either the $u$ or the $v$ directions. We start from the data form of the polynomial $\Delta_L(U,V,e,p)$ given in Example \ref{plainweavedata}, and in each case repeat the same variable $p$ for all the weft strands (in the $u$ direction), and $e$ for all the warp strands, (in the $v$ direction) to arrive at a formula for the multiple cell polynomial.

When doubled in the $v$ direction, using the cell determined by $u$ and $w=2v$ this gives the polynomial $Q(U,W,e,p)=\Delta_L(U,V,e,p)\Delta(U,-V,e,p)$ with $W=V^2$.

Here $$\Delta_L=-(1-U)^2(1+V^2)+[2+(1-e)^2(1-p)^2U/ep+2U^2]V,$$ giving
$$Q=(1-U)^4(1+W)^2-[2+(1-e)^2(1-p)^2U/ep+2U^2]^2W.$$

The Alexander data can be  presented  on a
 translation lattice  on the fabric, after expansion in terms of $U$ and $W$. The Salkeld factorisation characteristic of the multiple cell repeat is not immediately obvious in the lattice, but shows up from the Alexander polynomial when the variable $W$ is replaced by $V^2$.

\end{example}

In fact the plain weave itself as given traditionally is invariant under a larger group of translations than those generated by $u$ and $v$. The translations by $\ha(u+v)$ and $\ha(u-v)$ generate the full group of invariant translations, and give rise to half-size repeating cells. The kernel based on the choice of $u,v'=\ha(u+v)$ can be used to calculate the Alexander data for this larger group of translations. If we use the translations $\ha u$ and $\ha v$ to index rows and columns of a matrix then the Alexander data has the simple form shown, where as before we write $e=a^2$ and $p=b^2$.

\begin{figure}[ht!]
\bc
$ \begin{array}{ccc}
-1&.&1\\[5mm]
.&(a-a^{-1})(b-b^{-1})&.\\[5mm]
1&.&-1
\end{array}$
\ec \caption{Plain Weave minimal cell data} \label{plainweaveminimaldata}
\end{figure}

In this matrix view of the Alexander data the entries indicated by a dot (.) correspond to translations which do not belong to the group.
The Salkeld factorisation of the traditional polynomial arising from the fact that the repeat lattice is not minimal  can be seen from its data form in terms of $U$ and $V$  by putting $UV=V'^2$.

\begin{example} We can also see examples of the multiple cell repeat when we study twills.  Figure \ref{twillcellfig} shows a picture of a $1/2$ twill, with the traditional repeating cell based on vectors $u$ and $v$.
\begin{figure}[ht!]
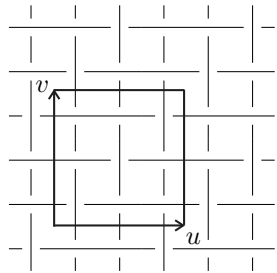

\bc
{\labellist
\pinlabel {$u$} at 178 32
\pinlabel {$v$} at 32 178
\endlabellist\twillcell} \ec
\caption{$1/2$ twill in traditional representation} \label{twillcellfig}
\end{figure}

The complete translation symmetry group is generated by $u$ and ${v'=\frac{1}{3}(u+v)}$, giving a repeat cell of one third the size,  shown in Figure \ref{twillmincellfig}.

\begin{figure}[ht!]
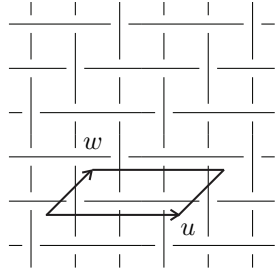

\bc
{\labellist
\pinlabel {$u$} at 173 37
\pinlabel {$w$} at 81 120
\endlabellist\twillmincell} \ec
\caption{Minimal unit cell of $1/2$ twill} \label{twillmincellfig}
\end{figure}

 The data for the kernel based on $u$ and $v'$ is given here, with $u$ drawn horizontally and $v'$ vertically.

\bc
$\begin{array}{ccc}
\sqrt{ab}&-\sqrt{b/a}&0\\[4mm]
0&\sqrt{b/a}(a^{-1}-a)(b^{-1}-b) &0\\[4mm]
0&\sqrt{a/b}(a^{-1}-a)(b^{-1}-b) &0\\[4mm]
0&-\sqrt{a/b}&1/{\sqrt{ab}}
\end{array}$
\ec

When  placed on the plane of the fabric  the data, drawn with $\frac{1}{3}u$ horizontally and
$\frac{1}{3}v$ vertically, appears as shown in Figure \ref{twillminimaldata}.

\begin{figure}[ht!]
\bc $\begin{array}{cccc}
\sqrt{ab}&.&.&-\sqrt{b/a}\\[4mm]
.&.&\sqrt{b/a}(a^{-1}-a)(b^{-1}-b)\\[4mm]
.&\sqrt{a/b}(a^{-1}-a)(b^{-1}-b)&.&.\\[4mm]
-\sqrt{a/b}&.&.&1/{\sqrt{ab}}
\end{array}$
\ec \caption{$1/2$ twill minimal cell data} \label{twillminimaldata}
\end{figure}

In this array a dot (.) is used for entries which do not correspond to any translation symmetry, as in the minimal cell display for plain weave shown in Figure \ref{plainweaveminimaldata}.

As in the case of the minimal cell repeat for plain weave, discussed above, Salkeld's theorem shows that the data form of the Alexander polynomial for the kernel based on $u$ and $v$ will factorise when we write $UV=V'^3$, corresponding to the relation $u+v=3v'$. One of the resulting factors, when written in terms of $U$ and $V'$, gives the data form for the minimal cell based on $u$ and $v'$. A direct check starting from the Alexander polynomial of the kernel based on $u$ and $v$ does indeed confirm this.

\end{example}

More generally, an $m/n$ twill consists of a weave in which there are $m$ warp overlaps and $n$
weft overlaps on any warp or weft thread within the repeat, with an offset of one thread in each
successive row. Take $u$ and $v$ as horizontal and vertical translations by $m+n$ threads. These
are in the symmetry group of the fabric, but the complete symmetry group can be generated by $u$
and ${v'=\frac{1}{m+n}(u+v)}$, giving a minimal repeat cell of $\frac{1}{m+n}$ times the obvious
cell generated by $u$ and $v$. Plain weave itself can be regarded in this context as a $1/1$ twill.

\subsection{Layered fabrics}\label{layered}

If a fabric decomposes into two separate layers then the plane separating the two layers
becomes a torus separating a kernel $N$ for the fabric into some curves $W_1,\ldots,W_r$
on one side of the thickened torus, coming from the lower layer of the fabric, and
others $T_1,\ldots,T_k$ on the other side coming from the top layer. The face curve $X$
lies on the side of the torus  containing the curves $T$, while the back curve $Y$ lies
on the side containing the curves $W$.  Because of the torus separating the components
of $N$ there is a decomposition of the polynomial $\Delta_N$ as a product. This arises
as a special case of a general result of Fox, best described in the following `Fox gluing
formula'.

Suppose that we have two links $L=T_1\cup\ldots\cup T_k\cup A$ and
$M=W_1\cup\ldots\cup W_r\cup B$. Remove a neighbourhood $V_A$ of $A$ and $V_B$ of $B$
and glue $S^3-V_A$ to  $S^3-V_B$, matching the meridian $a$ of $V_A$ to the longitude
of $V_B$ and the meridian $b$ of $V_B$ to the longitude of $V_A$.  The curves
$T_1\cup\ldots\cup T_k\cup  W_1\cup\ldots\cup W_r$ form a link $N$ in the resulting
manifold, which is again $S^3$ if one of $A$ or $B$ is unknotted.

\begin{theorem}[Fox gluing formula] (\cite{FoxV}) \label{FG} The Alexander polynomial of
the link $N$ resulting from this gluing is given by
$$\Delta_{N}=\Delta_L \Delta_M$$ after substituting $a=<B>$ in ${\bf Z}[w_j^{\pm1}]$
and $b=<A>$ in ${\bf Z}[t_i^{\pm1}]$.
\end{theorem}
\begin{remark} Theorem \ref{TF} is a corollary of this formula. Take the link $M$ to be the
unknot $B$ with no further components.  Then $<B>=1$ and $\Delta_M=\ds\frac{1}{1-b}$. The link
$N$ is the link $L-A$, and Theorem \ref{FG} shows that $(1-<A>)\Delta_{L-A}=\Delta_L|_{a=1}$.
\end{remark}

\begin{theorem}\label{layers}
If a  fabric kernel $N$ arises from a layered fabric as above, then its Alexander polynomial factorises as
$$\Delta_N=\Delta_L(x_L,y_M\prod w_j^{a_j},t_1,\ldots,t_k)\Delta_M(x_L\prod t_i^{b_i}, y_M,w_1,\ldots, w_r),$$ where the numbers $a_j$ are the linking numbers of the fabric strands in the bottom layer with the face curve, while the numbers $b_j$ are the linking numbers of the strands in the top layer with the back curve.
\end{theorem}

\begin{proof}
If a  fabric kernel $N$ arises from a layered fabric
then $N$  comes from  a gluing
construction to which Theorem \ref{FG} applies. Define two other   links which correspond to
the two layers of the fabric. These are the kernel for the top layer, $L=X_L\cup Y_L\cup T_1\cup\ldots\cup T_k$, given by
deleting the curves $W$ from $N$, and the kernel for the bottom layer, $M=X_M\cup Y_M\cup W_1\cup\ldots\cup W_r$, given by
deleting the curves $T$.  Then the link $N$ comes from $L$ and $M$ by the gluing construction,
taking $A=Y_L$ and $B=X_M$.

Now $<Y_L>=x_L\prod t_i^{b_i}$ and $<X_M>=y_M\prod w_j^{a_j}$. In $N$ the auxiliary curves
are $X=X_L$, the face curve of the top layer, and $Y=Y_M$, the back curve of the bottom layer.  Hence the Alexander polynomial factorises as
$$\Delta_N=\Delta_L(x_L,y_M\prod w_j^{a_j},t_1,\ldots,t_k)\Delta_M(x_L\prod t_i^{b_i}, y_M,w_1,\ldots, w_r).$$

Assuming that the fabric has no closed components this factorisation is non-trivial, since
evaluation of either factor when $t_i=w_j=1$ for all $i$ and $j$ gives a polynomial which
is not $0$ or $1$.
\end{proof}

\begin{corollary}
If a fabric decomposes into $r$ layers then the polynomial of
any kernel factorises into $r$ factors corresponding to the kernels of the layers.
\end{corollary}
The number
of non-trivial factors in the Alexander polynomial of any kernel is then an upper bound for the number of layers.
It follows that if the polynomial $\Delta_N$ does not factorise then the fabric does not
decompose into layers.

The converse does not hold, since some fabrics which do not decompose
into layers can have kernels whose polynomials do factorise non-trivially.

\subsection{Strip-like fabrics}

A fabric  may decompose into disconnected parallel strips. When this happens there will be only
one axial type for the strands. Assume that one of the generators $u$ is chosen in this direction.
 In this case the polynomial for the kernel will be independent of the variable
$x$ and must have a factor of $1-y$ even before setting all the fabric variables to $1$.  It is
worth testing a fabric with only one axial type for the possibility that it may fall into
strips. While a knitted fabric will only have one axial type it may be possible to confirm
that it is a genuinely 2-dimensional fabric if its polynomial has no factor of $1-y$, or has a
non-trivial dependence on $x$.  The knitted fabrics studied in examples \ref{jerseyexample}, \ref{jerseybackexample}, \ref{jerseydoublingexample}, \ref{ribexample} and \ref{garterexample}, with $u$ as their axial type, do indeed exhibit a factor $1-y$, but the polynomials depend non-trivially on $x$, confirming that the fabrics do not fall into strips.

\begin{figure}[ht!]
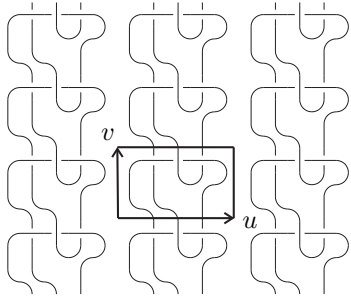

\bc
{\labellist
\pinlabel {$u$} at 327 332
\pinlabel {$v$} at 190 417
\endlabellist\chaincell} \ec
\caption{Warp-knitted chain} \label{chainfig}
\end{figure}

The warp-knitted chain illustrated in Figure \ref{chainfig} shows a linear chain repeated so
as to give a doubly periodic fabric, in this case presented with a single axial type $v$ .  The Alexander polynomial of its kernel,  $$
(x-1)(xt+1-t)(xt-x+1)
,$$ depends on $x$ and $t$ only, and has
a factor of $1-x$, both of which are pointers to the strip-like nature of the fabric. Its data form then depends on $V$ only and the data appears as a linear array of polynomials in the $v$ direction.

\begin{figure}[ht!]
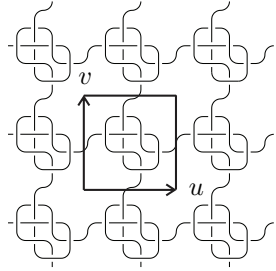

\bc
{\labellist
\pinlabel {$u$} at 182 75
\pinlabel {$v$} at 75 182
\endlabellist\fishnetcell} \ec
\caption{Fishing net} \label{fishnetfig}
\end{figure}

An example of a fabric which is not a traditional knit or weave is the fishing net, shown in
Figure \ref{fishnetfig}.  It has one axial type,
$u+v$, relative to the choice of generating translations indicated, but it does not have the
geometric characteristics of a knitted fabric with this as the thread direction.

The Alexander polynomial of the fishing net kernel based on this choice of generators is
$$
t^2x-t^3x+t^2y-t^3y-4t^3xy+2t^2xy-txy-t^5xy+3t^4xy+t^4-2t^3+4t^2-3t+1.
$$
To convert to data form we set $U=ys, V=xs$, with $t=s^2$. This leads to the display of Alexander data,
\bc
$\begin{array}{cc}
s^{-1}-s &-s^4+3s^2-4+2s^{-2}-s^{-4}\\
s^4-2s^2+4-3s^{-2}+s^{-4} & s^{-1}-s,
\end{array}$
\ec
confirming that this is not a strip-like fabric.

\subsection{Vassiliev invariants}
In a recent paper, \cite{GV},  Grishanov, Meshkov and Vassiliev have looked at the use of Vassiliev invariants of a curve in the thickened torus or Klein bottle as a tool for distinguishing textile patterns.

Following the paper of \cite{morton} in relation to Fiedler's invariant we believe that the Alexander polynomial of  a kernel can be used to find some of the Vassiliev invariants for a fabric with a single fabric curve.

Suppose that a kernel $L$ has been constructed from the fabric, with polynomial $\Delta_L(x,y,t)$, normalised to have Torres symmetry. Expand $\Delta_L(x,y,e^h)$ as a power series $\sum a_r(x,y) h^r$.  The coefficients can be regarded as elements of the integer group ring of $H_1(T^2)$, where $x$ and $y$ are represented by the curves $V$ and $U$ respectively in the torus $T^2$. In this setting the polynomial $a_r$ appears to be a Vassiliev invariant of degree $r$ in the sense of \cite{GV}. The constant term $a_0$ is, by Theorem \ref{axial}, the polynomial $1-x^a y^b$ up to normalisation, and the term $x^ay^b$ represents the homotopy class of the fabric component. This is constant over all homotopic curves in the torus, as for the degree $0$ invariant of  \cite{GV}.

The invariants of the fabric curve of degree  $1$ should also arise in this way.
A check on  the examples above with only one fabric component give the following results for $a_1$.

\begin{tabular}{lll}
Chain mail&
$a_0=0$ &$a_1= x+x^{-1}-y-y^{-1}$\\
Single jersey (face) &
$a_0=1-y$ &$a_1=a_0(x^{-1}-x)$\\
Single jersey (back)&
$a_0=1-y$ & $a_1=a_0(x-x^{-1})$\\
$m \times n$ rib &
$a_0=1-y$&  $a_1=a_0(n-m)(x-x^{-1})$\\
Fish net&
$a_0=1-xy $ & $a_1=\ha-x-y+\ha xy$\\
Warp knitted chain&
$a_0=1-x$&$a_1=a_0(x+1-x^{-1})$
\end{tabular}
\section{Further investigations} Most of our calculations have made use of a Maple program based on \cite{Morton} for calculating the multi-variable Alexander polynomial of a closed braid. This has required work to present a kernel for the fabric in closed braid form.  We hope to develop a similar program to handle links presented as plats, or even just from the diagram of the fabric inside a unit cell. This approach looks likely to relate to  a presentation of the diagram on the torus as a genus $1$ virtual knot, \cite{Kauffman}.

\section*{Acknowledgments}
We would like to thank De Montfort University and the Mathematics and Modelling Research Centre at the University of Liverpool for support during this work.

\end{document}